\documentclass[12pt]{article}

\usepackage[usenames]{color}

\def\chi{{\mathcal{X}}}

\def\one{\mbox{${\bf1}$}}

\def\calx{{\mathcal{X}}}

\def\cala{\mathcal{A}}
\def\({\left(}
\def\){\right)}
\def\pf{\n{\bf Proof.} }
\def\vsp{\vspace*{1,5mm}\\ }
\def\vspp{\vspace*{2mm}\\ }

\def\bk{\bigskip }
\def\mk{\medskip }
\def\sk{\smallskip }

\def\n{\noindent }
\def\dd{\displaystyle}
\def\hf{{\hfill\rule{1,5mm}{1,5mm}}}
\def\D{{\Delta}}

\def\barr{\begin{array}}
\def\earr{\end{array}}

\def\bit{\begin{itemize}}
\def\itemi{\item[{\rm(i)}]}
\def\itemii{\item[{\rm(ii)}]}
\def\itemiii{\item[{\rm(iii)}]}

\def\eit{\end{itemize}}

\def\D{{\Delta}}
\def\FP{Fokker--Planck}

\usepackage{amsmath, amsthm, amscd, amsfonts, amssymb}

\newtheorem{theorem}{Theorem}[section]

\newtheorem{corollary}[theorem]{Corollary}

\newtheorem{lemma}[theorem]{Lemma}
\theoremstyle{definition}

\newtheorem{remark}[theorem]{Remark}

\def\Lbb{{\Lambda}}
\def\1{^{-1}}
\def\vsp{\vspace*{2mm}\\ }

\def\calf{{\mathcal{F}}}

\def\calx{{\mathcal{X}}}

\def\rr{{\mathbb{R}}}
\def\nn{{\mathbb{N}}}
\def\9{{\infty}}

\def\lbb{{\lambda}}
\def\a{{\alpha}}
\def\b{{\beta}}

\def\wt{\widetilde}
\def\ov{\overline}
\def\vf{{\varphi}}
\def\oo{{\omega}}
\def\ooo{{\Omega}}
\def\pp{{\partial}}
\def\D{{\Delta}}
\def\vp{{\varepsilon}}
\def\barr{\begin{array}}
\def\earr{\end{array}}

\def\dd{\displaystyle}
\def\bk{\bigskip }
\def\sk{\smallskip}
\def\n{\noindent }

\def\pas{\mathbb{P}\mbox{-a.s.}}

\def\vsp{\vspace*{2mm}\\ }
\def\ff{\forall }
\def\({\left(}
\def\){\right)}
\def\<{\left<}
\def\>{\right>}

\title{Nonlinear  Fokker--Planck equations driven~by~Gaussian linear multiplicative noise}
\author{Viorel Barbu\thanks{Octav Mayer Institute of Mathematics of  Romanian Academy,     Ia\c si, Romania.  Email: vbarbu41@gmail.com}\and Michael R\"ockner\thanks{Fakult\"at f\"ur Mathematik, Universit\"at Bielefeld,  D-33501 Bielefeld, Germany.  Email: roeckner@math.uni-bielefeld.de}}
\date{}

\begin{document}
\maketitle
\begin{abstract}
\n Existence   of a strong solution in $H^{-1}(\rr^d)$  is proved for  the stochastic nonlinear Fokker--Planck equation  $$dX-{\rm div}(DX)dt-\Delta\beta(X)dt=X\,dW \mbox{ in }(0,T)\times\rr^d,\ X(0)=x,$$respectively, for a corresponding random differential equation. Here $d\ge1$,  $W$ is a  Wiener process  in $H^{-1}(\rr^d)$, $D\in C^1(\rr^d,\rr^d)$ and $\beta$ is a   continuous monotonically increasing function satisfying some appropriate sublinear growth conditions  which are compatible with the phy\-si\-cal models arising  in statistical mechanics. The solution exists for $x\in L^1\cap L^\9$ and preserves positivity.  If $\b$ is locally Lipschitz, the solution is unique, pathwise  Lipschitz  continuous with respect to initial data in $H^{-1}(\rr^d)$.  Stochastic {Fokker-Planck} equations with non\-linear drift of the form $dX-{\rm div}(a(X))dt-\Delta\beta(X)dt=X\,dW$ are also considered for   Lipschitzian {continuous} functions $a:\rr\to\rr^d$.\sk
\end{abstract}

\section{Introduction}\label{s1}
We  first consider the stochastic partial differential equation
\begin{equation}\label{e1.1}
\barr{l}
dX-{\rm div}(DX)dt-\Delta\beta(X)dt=X\,dW \mbox{ in }(0,T)\times\rr^d,\ T>0,\\
X(0,\xi)=x(\xi),\ \xi\in\rr^d,\ 1\le d<\9,
\earr\end{equation}
where $W$ is a   Wiener process in \mbox{$H^{-1}:=H^{-1}(\rr^d)$} over a stochastic basis $(\ooo,\calf,(\calf_t)_{t\ge0},\mathbb{P})$ with normal filtration $(\calf_t)_{t\ge0}$ of the form
\begin{equation}\label{e1.2}
W=\sum^N_{j=1}\mu_je_j\beta_j.\end{equation}
Here $\{e_1,...,e_N\}$ is an orthonormal system in $H^{-1}(\rr^d)$ belonging to $ C^2_b(\rr^d)\cap W^{2,1}(\rr^d)$, $\mu_j\in\rr$ and $\{\beta_j\}^\9_{j=1}$ are independent $(\mathcal{F}_t)$-Brownian motions on $(\ooo,\mathcal{F},\mathbb{P})$. As regards the functions $D:\rr^d\to\rr^d$ and $\beta:\rr\to\rr$,   we assume that
{\it\bit\itemi $D\in C^1_b(\rr^d;\rr^d);\ |D|\in L^1(\rr^d),\ {\rm div}\,D\in  L^2(\rr^d).$
\itemii $\beta\in  C(\rr)\cap C^2(\rr\setminus\{0\})$ is
     monotonically nondecreasing,   $\beta(0)=0$, and there are $m\in[0,1]$, $a_i\in(0,\9)$, $i=1,2,3$, such~that
    \begin{eqnarray}
    |\beta(r)|&\le&a_1|r|^m,\ \ff r\in\rr,\label{e1.3}\\
    |\beta''(r) r^2| +
    \beta'(r)|r|&\le&a_2|\beta(r)|,\
    \ff r\in\rr\setminus\{0\},\label{e1.4}\\
    \beta'(r)&\ne&0\mbox{ and }{\rm sign}\,r\,
    \beta''(r)\le0,\ \ff r\in\rr\setminus\{0\}.
    \quad\label{e1.5}
        \end{eqnarray}
        \itemiii There exists  a decreasing function $\vf:(0,1]\to(0,\9)$ such that
        \begin{equation}\label{e1.7}
        \beta'(\lbb r)\le\vf (\lbb)\beta'(r),\ \ff r\in\rr\setminus\{0\},\ \lbb\in (0,1].\end{equation}
               \eit}
       We note here that since, by \eqref{e1.5},  $\beta'$ is decreasing on $(0,\9)$ and increasing on $(-\9,0)$, we also have
       \begin{equation}\label{e1.7prim}
       \beta'(r)\le\beta'(\lbb r),\ \ff r\in\rr\setminus\{0\},\ \lbb\in(0,1].
       \end{equation}

     A typical example is $\b(r)\equiv a_1 r|r|^{m-1},$ where $a_1>0.$

     It should be said that   $e^{\pm W}$ is a linear multiplier in the spaces $L^p$ and $H^{1}$ and this fact will be frequently used in the sequel.

         Equation \eqref{e1.1}, which in the linear, deterministic case (that is, for\break  $\beta(r) \equiv  ar,$  $W=0$) reduces to the classical Fokker--Planck equation, \mbox{describes} the par\-ticle transport dynamics in disordered media driven by highly irre\-gu\-lar or stochastic field forces. This is the so called anomalous diffusion dy\-na\-mics (see, e.g., \cite{9}, \cite{10}) in contrast to the normal diffusion processes governed by the linear \FP\ equation.
         
         %%%%%%%%%%%%%%%%%
         The stochastic version \eqref{e1.1} considered here can be viewed as a Fokker-Planck equation in a random environment or a generalized mean field Fokker-Planck equation (\cite{9ax}, \cite{10a}, \cite{10aa}).  
         %%%%%%%%%%%%%%%

    The case considered here, that is hypothesis \eqref{e1.3} with $0\le m\le1$ is that of a fast diffusion (see, e.g., \cite{5}) which, for  $D\equiv0$  is relevant in plasma physics and the kinetic theory of gas. It should be said that in   statistical physics, the deterministic \FP\ equation \eqref{e1.1} is related to the so-called correspondence principle (see, e.g., \cite{10}, \cite{11}) in statistical mechanics which associates this equation to the entropy function
    $$S(u)=\int_{\rr}\Phi(u)d\xi,$$ where the function $\Phi\in C(\rr)\cap C^2(\rr\setminus\{0\})$ satisfies
        \begin{equation}\label{e1.10}
    \Phi''<0,\ \Phi'\ge0,\ \Phi'(0)=+\9,\end{equation}
    and $\beta$ is defined by
    \begin{equation}\label{e1.11}
    \beta(r)=\Phi(r)-r\Phi'(r),\ \ff r\ge0.\end{equation}
        For instance, if $\beta(r)\equiv a\ {\rm sign}\ (r)\log(1+|r|)$, $a>0$, and $\Phi(u)=-u\log u+$ $(1+u)\log(1+u)$, then \eqref{e1.1} is the classical boson  equation in the Bose--Einstein statistics (see, e.g., \cite{10}), while for $\beta(r)\equiv a|r|^{m-1}r$, one gets the so-called Plastino and Plastino model \cite{11} in statistical mechanics.

        We note that  in  both cases  $\b$ satisfies (ii) and (iii) above, and in the first case $\b$ is locally Lipschitz.  

        Assumption (ii) leaves out the low diffusion case $m>1$ which is relevant in  porous media dynamics of low diffusion processes. (See, e.g., \cite{5}.)  However, for the examples in statistical mechanics mentioned above, the case $m>1$ is not relevant. In fact, the entropy function corresponding to $\beta(u)=u^m$ is by \eqref{e1.11} formally given in $1-D$ by
        $$S(u)=\frac1{1-m}\int_{\rr}(u^m-u)d\xi,\ \Phi(u)=\frac1{1-m}\,(u^m-u),$$
        for which   the entropic conditions \eqref{e1.10} are not satisfied if $m>1.$

    For   vanishing drift $D$, equation \eqref{e1.1} reduces to the fast diffusion stochastic porous media equation studied in \cite{7a} (see, also, \cite{5}).

    By the transformation
    \begin{equation}\label{e1.12}
    X(t)=e^{W(t)}y(t),\ t\ge0,\end{equation}
    equation \eqref{e1.1} reduces, via It\^o's formula, to the random differential equation (see, e.g.,  \cite{4}, \cite{6}, \cite{7})  
    \begin{equation}\label{e1.13}
    \barr{l}
    \dd\frac{\pp y}{\pp t}-e^{-W}{\rm div}(e^WDy)-e^{-W}\Delta\beta(e^Wy)+\dd\frac12\,\mu y=0\mbox{ in }(0,T)\times\rr^d,\vspp
    y(0,\xi)=x(\xi),\ \xi\in\rr^d,\earr\end{equation}
    where
    \begin{equation}\label{e1.14}
    \mu=\sum^N_{j=1}\mu^2_je^2_j.\end{equation}
Here, without loss of generality, we assume that $t\mapsto W(t)(\omega)\in H^{-1}$ is continuous for all $\omega\in\Omega$.

    The purpose of this work is to show that, under hypotheses (i)-(iii),  for every $\oo\in\ooo$, $1\le d<\9$, and $x$ in a suitable space, the Cauchy problem \eqref{e1.13} has  at least one    strong solution  which is unique if, in addition, $\b$ is locally Lipschitz on $
    \rr$.   By a strong solution to \eqref{e1.13} we mean an absolutely continuous function  $y:[0,T]\to H^{-1}(\rr^d)$ such that ${\rm div}(e^{W}Dy)(t)\in H^{-1}$, a.e. $t\in(0,T)$, and \eqref{e1.13} holds on $ (0,T).$ Of course, if $y$ is $(\calf_t)_{t\ge0}$-adapted (which we shall show), then $X=e^Wy$ is a strong solution to \eqref{e1.1}. A nice feature of the random differential equation \eqref{e1.13} and its version with a nonlinear function in its divergence part (see equation \eqref{e4.2} below) is that, though it is not of accretive type in any of the spaces $H^{-1}(\rr^d)$ or $L^1(\rr^d)$, which are naturally associated with nonlinear parabolic equations of this type, it turns out to be accessible by  the theory of nonlinear semigroups of contractions in $L^1(\rr^d)$, by a modification of the Crandall-Liggett discretization scheme for perturbed nonlinear accretive equations (see Appendix). 
    
    However, the general existence theory for the nonlinear accretive Cauchy problem in a Banach space is not directly applicable to equation \eqref{e1.1} because $W$ is not smooth. So, the first step was to approximate  $W$ by a family of smooth random functions $\{W_\vp(t)\}_{\vp>0}$ and   so equation  \eqref{e1.13} too by a family of nonlinear evolution equation with smooth time-dependent coefficients (see equation \eqref{e3.2} below). Afterward, one passes to the limit $\vp\to0$ in the corresponding equation by combining sharp $H\1$-energetic and $L^1$-techniques. This approach which will lead to existence of a strong solution $y$ to \eqref{e1.13} is one of the main novelty of this work.

    In \cite{4}, the authors studied equation \eqref{e1.13} for $m\in(1,5) $ and $1\le d\le 3$, on a bounded domain in the special case of a vanishing drift term $D$. It~should be said, however, that the treatment  in $\rr^d$ developed here is  quite  different and requires specific techniques to be made precise below. (Under related hypotheses on $\beta$, the existence for the stochastic equation \eqref{e1.1} with $D\equiv0$ was also studied in \cite{7a}.)

    In \cite{11aa},   the following  parabolic-hyperbolic quasilinear stochastic equation  was recently studied on $T^d$ in the framework of kinetic solutions
    \begin{equation}\label{e1.15}
    dX-{\rm div}(B(X))dt-{\rm div}(A(X)\nabla X)dt=\Phi(X)dW,\end{equation}
    where $B\in C^2(\rr,\rr^{d\times d})$ and $A\in C^1(\rr;\rr^{d\times d})$. (Along these lines, see also \cite{11aaa}.) It should be said, however, that there is no overlap with our work   as far as conditions (i) on the nonlinear diffusion term $\beta$ is concerned for which one assumes here different conditions to cover fast diffusions. In fact, the results of \cite{11aaa}, though obtained in a more general context, apply to low diffusion equations (that is, $\beta(r)\approx ar^m$, $m\ge2$, $a(r)\approx r^k$, $k>1$). In addition, the rescaling technique used here is different from  that used in \cite{11aaa} and its main advantage is that it leads to sharper regularity results for solutions by fully exploiting the parabolic nature of the resulting random differential equation.

    \section{Notation and the main results}
    \setcounter{equation}{0}

We shall denote  the norm of the space $\rr^d$ by $|\cdot|$ and by $\<\,,\,\>$ the Euclidean inner product.
    Let $L^p(\rr^d)=L^p$,  $1\le p\le\9,$ denote the standard real $L^p$ space on $\rr^d$ with  Lebesgue measure. The scalar product of $L^2$ is denoted by $(\cdot,\cdot)_2$. The norm of $L^p$ will be denoted by $|\cdot|_p$.  $H^1(\rr^d)$, briefly denoted $H^1$, is the Sobolev space $\left\{u\in L^2;\  \frac{\pp u}{\pp \xi_i}\in L^2,\ i{=}1,2,...,d\right\}$ with   the standard norm $\|u\|_{H^1}=  \left(\int_{\rr^d}(u^2+|\nabla u|^2) d\xi\right)^{\frac12}$. The dual space of $H^1$ will be denoted by $H^{-1}$ and its norm by $|\cdot|_{-1}$. Likewise, $W^{r,p}=W^{r,p}(\rr^d)$, $r\in\mathbb{N}$, $p\in[1,\9]$, denote the usual Sobolev spaces. Denote by $\Delta$ the Laplace operator on $\rr^d$. By~$W^{1,p}([0,T];H^{-1})$ we denote the space of all absolutely continuous\break $u:[0,T]\to H^{-1}$ such that $u,\frac{du}{dt}\in L^p(0,T;H^{-1})$. Given a Banach space $X$, let   $L^p(0,T;X)$  denote the space of $X$-valued Bochner $L^p$-integrable functions on $(0,T)$. By $C([0,T];X)$, we denote the space of continuous functions $u:[0,T]\to X$ and by $C^1([0,T];X)$ the cor\-res\-pon\-ding space of continuously differentiable functions.

We set
 $$D_0=\{x\in L^1\cap L^\9\cap H^1;\ \b(x)\in H^1,\ \D x\in L^1,\ \D\b(x)\in L^1\}.$$

\begin{lemma}\label{l20} Let $p\in[1,\9)$
and $x\in L^1\cap L^\9$.
Then there exist $u_n\in D_0$,
$n\in\nn$, such that $u_n\to x$
in $L^p$ and $\{x_n;\ n\in\nn\}$ is
bounded in $L^1\cap L^\9$. In particular,
$$\ov D^{L^p}_0=L^p,\quad
\ov D^{H^{-1}}_0=H^{-1},$$
where the left hand sides denote the closures of $D_0$ in the respective spaces.
\end{lemma}

\pf Because $L^2$ is dense in $H^{-1}$, it suffices to prove
 $$L^1\cap L^\9\subset\ov D^{L^p}_0.$$ 
So, let $x\in L^1\cap L^\9$ and define
\begin{equation}\label{e20}
u(\xi)=\vf(\xi)e^{-\delta|\xi|^2},\ \xi\in\rr^d,
\end{equation}
where $\vf\in C^2_b(\rr^d),
|\vf|\ge\vp,\ \vp,\delta\in(0,1).$
Then, by \eqref{e1.3},
 $\b(u)\in L^1\cap L^\9$ and
$$\nabla\b(u)=\frac1\vf\,\b'(u)u(\nabla\vf-2\delta\vf\xi),$$
which is in $L^1\cap L^\9$ by \eqref{e1.3}, \eqref{e1.4}. So, $\b(u)\in H^1$. Furthermore, obviously,
$\D u\in L^1\cap L^\9$, and
$$\barr{lcl}
\D\b(u)&=&\dd\frac1\vf\,\b'(u)u[\D\vf-(2d\delta-4\delta^2|\xi|^2)\vf
-4\delta\xi\cdot\nabla\vf]\vsp
&&+\dd\frac1{\vf^2}\,\b''(u)u|
\nabla\vf-2\delta\vf\xi|^2.\earr$$
Since $|\vf|\ge\vp$, it follows by \eqref{e1.3} and \eqref{e1.4} that
 $\D\b(u)\in L^1\cap L^\9.$ We have
$$x=\lim_{\delta\to0}\lim_{\vp\to\9}(x^+\vee\vp-x^-\wedge(-\vp))e^{-\delta|\xi|^2},$$where
 both limits are in $L^p$ and,
 obviously, each function on the right under the limits for fixed $\vp,\delta\in(0,1)$ can be approximated by functors of type \eqref{e20} in $L^p$.\hf\bk

    Theorem \ref{t2.1} is the main result.
\begin{theorem}\label{t2.1} Under Hypotheses {\rm(i)--(iii)}, for each $x\in  D_0$,   equation \eqref{e1.13} has, for each $\oo\in\ooo$,  at least one  strong  solution	
\begin{eqnarray}
&y\in W^{1,2}([0,T];H^{-1})\cap L^\9((0,T)\times\rr^d)\cap L^\9(0,T;L^1),\label{e2.1}\\[3pt]
&y\in   L^2(0,T;H^1),\label{e2.2}\\[3pt]
&    \beta(e^Wy)\in L^2(0,T;H^1).\label{e2.3}\end{eqnarray}
Moreover, if $x\ge0$, a.e. on $\rr^d$,
    then $y\ge0$, a.e. on $(0,T)\times\rr^d$.

    If   $\b$ is locally     Lipschitz on $\rr$ and assumptions {\rm(i)--(iii)} hold, then there is a unique  strong   solution $y$ to \eqref{e1.13}. This solution is $(\calf_t)$--adapted, the map $D_0\ni x\to y(t,x)$ is Lipschitz     from $H^{-1}$ to $C([0,T];H^{-1})$ on balls
    in $L^1\cap L^\9$ and   $y$ extends by density to a strong solution to \eqref{e1.13}, satisfying \eqref{e2.1}, \eqref{e2.3}, for all $x\in L^1\cap L^\9.$
    \end{theorem}

Now, coming back to equation \eqref{e1.1}, we recall (see, e.g., \cite{5}, \cite{4}, \cite{7a}) that a continuous $(\calf_t)_{t\ge0}$-adapted process $X:[0,T]\to H^{-1}$ is called strong solution to \eqref{e1.1} if the following conditions hold:
\begin{eqnarray}
 X\in L^2([0,T];L^2),&& \pas,\label{e2.4}\\[3pt]
\beta(X)\in L^2(0,T;H^1),&& \pas,\label{e2.5}\end{eqnarray}
\begin{equation}
\barr{r}
\dd X(t)-\int^t_0{\rm div}(DX(s))ds-\dd\int^t_0\Delta\beta(X(s))ds=x+\int^t_0X(s)dW(s),\vsp \ff t\in[0,T],\ \pas\earr \label{e2.6}
\end{equation}
We note here that, by \eqref{e2.4} and \eqref{e3.6prim} below,
$${\rm div}(DX)\in L^2(0,T,H^{-1}),\ \pas$$
The stochastic (It\^o-) integral in \eqref{e2.5} is the standard one (see \cite{8}, \cite{16a}, \cite{11a}). In fact, in the terminology of these references, $W$ is a $Q$-Wiener process $W^Q$ on $H^{-1}$, where $Q:H^{-1}\to H^{-1}$ is the symmetric trace class operator defined~by
$$Qh:=\sum^N_{k=1}\mu_k(e_k,h)_{-1}e_k,\ h\in H^{-1}.$$

\begin{theorem}\label{t2.2}  If $\b$ is locally Lipschitz on $\rr$ and  assumptions {\rm (i)--(iii)}  hold, then,
for every $x\in D_0$, equation \eqref{e1.1} has  a unique 
strong solution  $X=e^Wy$, which  satisfies
\begin{equation}\label{e2.6ax}
X e^{-W}\in W^{1,2}([0,T];H^{-1}),\ \pas,\end{equation}
and  $X\ge0$, a.e. on $(0,T)\times\rr^d\times\ooo$
if $x\ge0$, a.e. on $\rr^d$. Moreover,  the map $x\mapsto X(t,x)$ is $H^{-1}$-Lipschitz from balls in $L^1\cap L^\9$ to $C([0,T];H^{-1})$.
\end{theorem}

The argument used to show that $X$ is a
strong solution to \eqref{e1.1}
is standard up to a stopping time
argument and very similar to that
from the works \cite{6}, \cite{7}
and so it will be omitted.

 It should be said that assumptions of Theorem \ref{t2.2} (that is, (i)-(iii) and $\b$ locally Lipschitz) hold for the boson equation
$$dX-{\rm div}(DX)dt-\D(\log(1+|X|))dt=X\,dW$$
and  for other significant models in statistical mechanics. However, it leaves out the Plastino \& Plastino model \cite{8aa} for which all we can prove is the existence of a strong solution to the corresponding random differential equation \eqref{e1.13}.

A result as Theorem \ref{t2.2}
was previously proved in \cite{7a}
for equation \eqref{e1.1} in the
special case of   vanishing drift
$D$ by a direct approximation approach
to the stochastic equation \eqref{e1.1}.
The approach used here, based on the random
differential equation \eqref{e1.13},
is completely different and leads to
sharper results. Indeed, by \eqref{e2.1},
it follows that besides \eqref{e2.4}
the solution $X$ to \eqref{e1.1} satisfies also \eqref{e2.6ax},
 which is, of course, a new result.

It should be emphasized that the random differential
equation \eqref{e1.13}
has an interest in itself
as a model for particles
dynamics driven by random
transport and diffusion coefficients (see, e.g., \cite{9ax}).
In particular, the convergence of this solution to a stationary state or, more generally, the existence of a random attractor is a problem of utmost importance for its physical signi\-fi\-cance  related to the so-called Boltzmann $H$-theorem (see \cite{10}, \cite{12}).  We note here that, if our solution is unique for every fixed $\oo$, which is proved in this paper if $\b$ is locally Lipschitz, then, since it solves a deterministic PDE with random coefficients, it satisfies the strict cocycle property, so gives rise to a random dynamical system. This is the first and a fundamental ingredient to prove the existence of a random attractor. However, the uniqueness of solutions $y$ to \eqref{e1.13} under assumptions (i)-(iii) remains an open problem.

\section{Proof of Theorem \ref{t2.1}}
\setcounter{equation}{0}

Below we fix $\oo\in\ooo$, but do not express it in the notation.

Let $\beta^\vp_j\subset C^1([0,T];\rr)$, $1\le j\le N$,  be defined by $\b^\vp_j(t)=(\one_{[0,\9)}\b_j*\rho_\vp)(t)$, where $\rho_\vp(t)\equiv\frac1\vp\,\rho\(\frac t\vp\)$ is a standard mollifier with $\rho\in C^\9_0(\rr)$, $\rho\ge0$. We set
$$W_\vp(t,\xi)=\sum^N_{j=1}\mu_j e_j(\xi)\beta^\vp_j(t),\ t\ge0,\ \xi\in\rr^d.$$
Then we have for its time derivative

$$(W_\vp)_t\in C([0,T]\times\rr^d)$$
and $$W_\vp(t,\xi)\to W(t,\xi)\mbox{ uniformly in $(t,x)\in[0,T]\times\rr^d$}$$
as $\vp\to0.$

For each $\vp\in(0,1]$, consider the approximating equation of \eqref{e1.13}
    \begin{equation}\label{e3.1}
    \barr{l}
    \dd\frac{\pp y_\vp}{\pp t}-e^{-W_\vp}{\rm div}(e^{W_\vp}Dy_\vp)-e^{-W_\vp}\Delta(\beta(e^{W_\vp}y_\vp)+\vp e^{W_\vp}y_\vp)\vsp
    \hfill+\vp e^{-W_\vp}\beta(e^{W_\vp}y_\vp)+\dd\frac12\,\mu y_\vp=0\mbox{ in }(0,T)\times\rr^d,\vspp
    y_\vp(0,\xi)=x(\xi),\ \xi\in\rr^ d.\earr\end{equation}
    Setting $z_\vp=e^{W_\vp}y_\vp$, we get the equation
    \begin{equation}\label{e3.2}
    \barr{l}
    \dd\frac{\pp z_\vp}{\pp t}-\Delta(\beta(z_\vp)+\vp z_\vp)-{\rm div}(Dz_\vp)+\vp\beta(z_\vp)\vspp
    \hfill+\(\dd\frac12\,\mu-(W_\vp)_t\)z_\vp=0\mbox{ in } (0,T)\times\rr^d,\vspp
    z_\vp(0,\xi)=x(\xi),\ \xi\in\rr^d.\earr\end{equation}We have

    \begin{lemma}\label{l3.1} Assume that $x\in H^1$
    such that $\beta(x)\in H^1$. Then, for each $\vp\in(0,1]$,
     equation \eqref{e3.1} considered on $H^{-1}$
     has a unique strong solution $y_\vp$
     (see the Appendix) satisfying
     \begin{eqnarray}
    y_\vp&\in&W^{1,\9}([0,T];H^{-1})
    \cap L^\9(0,T;H^1).\label{e3.3}
    \end{eqnarray}
   Moreover, if $x\in D(A_1)$ with $D(A_1)$ defined as in the claim following \eqref{e3.13prim} below, then $y_\vp\in C([0,T];L^1)$ and $z_\vp=e^{W_\vp}y_\vp$, obtained as the limit of the finite difference scheme \eqref{e5.10prim}, is a mild solution to \eqref{e3.2} in the space $L^1$.\end{lemma}

    \pf It suffices to prove that equation \eqref{e3.2} has a unique solution
  \begin{equation}\label{e3.5}
  z_\vp\in W^{1,\9}([0,T];H^{-1})\cap L^\9(0,T;H^1),\end{equation} and $\beta(z_\vp):[0,T]\to H^1$ is right continuous.

  Let us first prove existence and uniqueness of a solution to \eqref{e3.2} considered as an equation on $H^{-1}$. Define the operator $A:D(A)\to H^{-1}$ by
    \begin{equation}\label{e3.8}
    Az=-\Delta(\beta(z)+\vp z)+\vp\beta(z)-{\rm div}(Dz)+\dd\frac\mu2\,z,\end{equation}
    with the domain $D(A)=\{z\in H^1:\beta(z)\in H^1\}$. We endow the space $H^{-1}$ with the scalar product
    $$\<y,z\>_{-1,\vp}={}_{H^1}\<(\vp I-\Delta)^{-1}y,z\>_{H^{-1}}
    %\int_{\rr^d}(\vp I-\Delta)^{-1}yz\,d\xi;
    \ y,z\in H^{-1},$$and with the corresponding norm $\|y\|_{-1,\vp}=(\<y,y\>_{-1,\vp})^{\frac12}.$ Taking into account that
    \begin{equation}\label{e3.6prim}
    \|{\rm div}(Dz)\|_{-1,\vp}\le\frac1{\sqrt{\vp}}
        \,|D|_\9|z|_2,\ \ff z\in L^2,\end{equation}
        we see that, for all $z,\bar z\in D(A)$,
    $$\<(A+\alpha I)z-(A+\alpha I)\bar z,z-\bar z\>_{-1,\vp}\ge0,$$
    if
    \begin{equation}\label{e3.8a}
    \alpha_\vp=\frac1\vp\,(|D|_\9 +\frac12\,|\mu|_\9).\end{equation}This means that $(A+\alpha I)$ is accretive in $H^{-1}$. Moreover, $A$ is quasi-$m$-accretive, that is, $R(\lbb+\alpha_\vp) I+A)=H^{-1}$ for all $\lbb>0$.   Indeed, for $f\in H^{-1}$, the equation
    \begin{equation}\label{e3.7a}
    (\alpha_\vp+\lbb)z-\Delta(\beta(z)+\vp z)+\vp\beta(z)-{\rm div}(Dz)+\frac\mu2\,z=f,\end{equation}
    or, equivalently,
    \begin{equation}\label{e3.9}
   \barr{l}
        \dd(\alpha_\vp{+}\lbb)(\vp I{-}\Delta)^{-1}z{+}\beta(z){+}\vp z{-}(\vp I{-}\Delta)^{-1}
        \({\rm div}(Dz)+\vp^2z-\frac\mu2\,z\)\vsp
        \qquad\qquad=(\vp I-\Delta)^{-1}f\earr \end{equation}
    has, for $\lbb>0$,   a unique solution $z\in L^2$.  Indeed, equation \eqref{e3.9} is of the form
    $$\vp z+B(z)+\Gamma z=(\vp I-\Delta)^{-1}f\in H^1,$$ where the operators $B:L^2\to L^2$ and $\Gamma:L^2\to L^2$ are given by
    $$\!\!\barr{rcl}
    B(z)(\xi)\!\!&{=}&\!\!\beta(z(\xi)),\mbox{ a.e. in }\rr^d,\vspp
    \Gamma(z)\!\!&{=}&\!\!\dd(\alpha_\vp+\lbb)(\vp I-\Delta)^{-1}z-(\vp I-\Delta)^{-1}\({\rm div}(Dz)+\vp^2z-\frac\mu2\,z\)\!.\earr$$
    Since $B$ is $m$-accretive and $\Gamma$ is accretive and continuous in $L^2$, it follows that $R(\vp I+B+\Gamma)=L^2$ and so there is a solution $z\in L^2$ to \eqref{e3.9}. Since, by \eqref{e3.9}, $\beta(z)+\vp z\in H^{1}$, since the inverse of $r\mapsto \beta(r)+\vp r$ is Lipschitz and equal to zero at $r=0$, it follows that $z\in D(A)$, as claimed.

    Now, we shall apply Lemma \ref{l5.1} and Corollary \ref{c5.2} in the Appendix, where $X=H^{-1},$ $A$ is the operator \eqref{e3.8} and $\Lambda(t)\in L(H^{-1}, H^{-1})$, $\ff t\in[0,T]$ defined~by
    \begin{equation}\label{e3.10}
    \Lambda(t)u=-(W_\vp)_t u,\ \ff u\in H^{-1},\end{equation}
    and get a strong solution $z_\vp$ to \eqref{e3.2} satisfying
    \begin{equation}\label{e3.10prima}
    z_\vp\in W^{1,\9}([0,T];H^{-1}).\end{equation}
    But, indeed, also
    $$z_\vp\in L^\9(0,T;H^1),$$
    i.e., \eqref{e3.5} holds. This can be seen as follows.

    By Corollary \ref{c5.2}, it immediately follows that
    \begin{equation}\label{e3.10secunda}
    \beta(z_\vp)+\vp z_\vp-(\vp z_\vp-\D)^{-1}{\rm div}(Dz_\vp)\in L^\9(0,T;H^1).\end{equation}
    An elementary consideration shows that, for $\vp\in(0,1)$,
     \begin{equation}\label{e3.10terta}
    |(\vp I-\D)^{-1}{\rm div}(Dz)|_{L^2}\le c|z|_{-1,\vp},\ \ff z\in L^2, \end{equation}
where $c$ is a constant (only depending on   $|D|_{C^1_{b}}$ and $d$). Since $z_\vp$ is a strong solution, we have $z_\vp\in D(A)\subset H^1(\subset L^2)dt$-a.e. Hence, it follows by \eqref{e3.10prima}-\eqref{e3.10terta} that
$$\beta(z_\vp)+\vp z_\vp\in L^\9(0,T;L^2),$$
hence also $z_\vp\in L^\9(0,T;L^2)$. So by \eqref{e3.6prim} we conclude
$$(\vp I-\D)^{-1}{\rm div}(Dz_\vp)\in L^\9(0,T;H^1).$$
Hence \eqref{e3.10secunda} implies that $\beta(z_\vp)+\vp z_\vp\in L^\9(0,T;H^1)$ and thus $z_\vp\in L^\9(0,T;H^1)$.

We are now going to construct  the realization   of the operator $A$ in $L^1.$
We consider the operator $A_0$ defined by
\begin{equation}
\label{e3.13prim}
\barr{c}
A_0z=-\D(\beta(z)+\vp z)+\vp\beta(z)-{\rm div}(Dz)+\frac\mu2\,z,\vsp z\in D(A_0)=D(A)\cap \{z\in L^1;\,\beta(z),\D(\beta(z)+\vp z)\in L^1\}.\earr
\end{equation}

\n{\bf Claim.} {\it Its closure  $A_1=\ov A_0$  in $L^1\times L^1$ is quasi $m$-accretive.}\mk

\n Indeed, since ${\rm div}\,D\in L^\9$, $D\in L^1\cap L^\9\subset L^2$, we have for all $z\in H^1\cap L^1$
 \begin{equation}\label{e3.10prim}
 \int_{\rr^d}{\rm div}(Dz){\rm sign}\,z\,d\xi=\int_{\rr^d}{\rm div}\,D|z|d\xi+\int D\cdot\nabla|z|d\xi=0.\end{equation}
 But, by \cite{1}, Theorem 3.5, also $D(A_0)\ni z\mapsto \D(\beta(z)+\vp z)$ is accretive on $L^1$; hence, since $\b$ is accretive, $A_0$ is  accretive on $L^1$ and hence so is $\ov A_0$. But we also have, for $\a>\a_\vp$,
 \begin{equation}\label{e3.10secund}
R(\a I+A_0)\supset H^{-1}\cap L^1,\end{equation}
because, for $f\in H^{-1}\cap L^1$, as we have seen above, there exists $z\in D(A)$ such that $\a z+Az=f$. But, indeed, $z\in L^1$. This can be seen as follows: for $\delta>0$, define for $r\in\rr$
\begin{equation}\label{e3.16az}
\chi_\delta(r):=\left\{\barr{rll}
1&\mbox{ if }&r>\delta,\vsp
\dd\frac r\delta&\mbox{ if }&r\in[-\delta,\delta],\vsp
-1&\mbox{ if }&r<\delta.\earr\right.\end{equation}
Then $\chi_\delta(z)\in H^1$ and, applying $_{H^1}\!\<\chi_\delta(z),\cdot\>_{H^{-1}}$ to \eqref{e3.7a}, we find
$$\barr{l}
\dd\a\int_{\rr^d}\chi_\delta(z)z\,d\xi
+\int_{\rr^d}\chi'_\delta(z)|\nabla z|^2(\b'(z)+\vp)d\xi\vsp
 \qquad\dd+\vp\int_{\rr^d}\chi_\delta(z)\b(z)d\xi
-\dd\int_{\rr^d}{\rm div}\,D\ \chi_\delta(z)z\,d\xi\vsp
 \qquad-\dd\int_{\rr^d}\<D,\nabla z\>\chi_\delta(z)d\xi
+\dd\frac12\int_{\rr^d}\chi_\delta(z)\mu\,z\,d\xi
=\dd\int_{\rr^d}\chi_\delta(z)f\,d\xi.\earr$$
Hence, dropping the second, third and sixth term (which are nonnegative) on the left hand side and then letting $\delta\to0$, because $D$, ${\rm div}\,D\in L^2$ we obtain
$$\a|z|_1\le|f|_1.$$
But then it   follows from \eqref{e3.9} that $\beta(z)\in L^1$ and hence, by \eqref{e3.7a}, that $z\in D(A_0)$ and \eqref{e3.10secund} is proved. Taking $L^1$-closure, we conclude that
$$\ov{R(\a I+A_0)}^{\,L^1}=L^1.$$
This implies that $\ov A_0$ is quasi-$m$-accretive, because for $\a$ large enough
$$R(\a I+\ov A_0)\supset\ov{R(\a I+A_0)}^{\,L^1},$$
  and the claim is proved.

 Then, again by Lemma \ref{l5.1} and Corollary \ref{c5.2}, applied to $X=L^1$ and to the operator $A_1$, it follows that for $x\in L^1$ equation \eqref{e3.2} has a unique mild solution $\wt z_\vp\in C([0,T];L^1)$ and $\wt y_\vp=e^{-W_\vp}\wt z_\vp$ is the mild solution to \eqref{e3.1}.

Let us note that $\wt z_\vp=z_\vp$ (and  $\wt y_\vp=y_\vp$, respectively) for $x\in D(A_0)$. Indeed, as seen in Lemma \ref{l5.1}, both $z_\vp$ and $\wt z_\vp$ are limits of finite difference scheme as \eqref{5.10}, where $A$ is given by \eqref{e3.8} and by $A_1=\ov A_0^{L^1}$, respectively. But, by \eqref{e3.10secund},
$$(I+hA_0)^{-1}y=(I+hA)^{-1}y,\ \ff y\in H^{-1}\cap L^1,\ \ff h\in(0,\alpha^{-1}_\vp).$$
The solutions $u_1\in L^1$ and $u\in H^1$ respectively of
\begin{equation}\label{e3.15prim}
u_1+h(A_1+\Lambda(ih))u_1=y\end{equation}
and
\begin{equation}\label{e3.15secund}
u+h(A+\Lambda(ih))u=y\end{equation}
for small enough $h$ are obtained by iterating the strict contractions\break $B_1:L^1\to L^1$, $B:H^{-1}\to H^{-1}$, defined by
$$B_1v:=(1+hA_1)^{-1}(y-h\Lambda(ih)v),\ v\in L^1,$$
and
$$Bv:=(1+hA)^{-1}(y-h\Lambda(ih)v),\ v\in H^{-1}.$$
Here $\Lambda(t)$ is given by \eqref{e3.10}, hence
$\Lambda(ih)$ leaves both $L^1$ and $H^{-1}$ invariant. Therefore, starting the iteration in a point $v_0\in H^{-1}\cap L^1$, we obtain by \eqref{e3.10secund} that
$$B^n_1v_0=B^nv_0\in D(A_0),\ \ff n\in\nn,$$
and that this sequence converges both in $L^1$ and $H^{-1}$. 

This  implies that
$$(I+h(A_0+\Lambda(ih)))^{-1}y=
(I+h(A+\Lambda(ih)))^{-1}y,\
i=0,1,...,\ \ff y\in H^{-1}\cap L^1.$$
This means that the finite difference schemes \eqref{e5.10prim} in Lemma \ref{l5.1},
 applied separately in the spaces $L^1$ and    $H^{-1}$,      lead for $x\in D(A)\cap D(A_1)$ to the same values $u^h=z^h_\vp$ ($\wt u^h=\wt z^h_\vp,$ respectively) and so, for the limit  $h\to0$, we get $z_\vp=\wt z_\vp$ for initial data $x\in D(A)\cap D(A_1)$. Hence $y_\vp=\wt y_\vp$, if $x\in D(A)\cap D(A_1).$~\hf\bk

To get rigorous estimates for solutions $y_\vp$ to equation \eqref{e3.1}, it is convenient to approximate it by the solution $y^\lbb_\vp$ to the equation
\begin{equation}\label{e3.11ab}
\barr{l}
\dd\frac{\pp y^\lbb_\vp}{\pp t}-e^{-W_\vp}{\rm div}(e^{W_\vp}Dy^\lbb_\vp)-e^{-W_\vp}\Delta(\beta_\lbb(e^{W_\vp}y^\lbb_\vp)\vspp
\qquad+\vp e^{W_\vp}y^\lbb_\vp)+\vp e^{-W_\vp}\beta_\lbb(e^{W_\vp}y^\lbb_\vp)+\dd\frac12\,\mu y^\lbb_\vp=0,\vspp
y^\lbb_\vp(0)=x,\earr\end{equation}
where $\beta_\lbb=\beta((I+\lbb\beta)^{-1})=\frac1\lbb\,(I-(I+\lbb\beta)^{-1})$ is the Yosida approximation of $\beta$. We recall that $\beta_\lbb$ is monotonically increasing, Lipschitzian and $$\lim_{\lbb\to0}\beta_\lbb(r)=\beta(r)\mbox{ uniformly on compacts in }\rr.$$We have

\begin{lemma}\label{l3.3a} For $\lbb\to0$, we have, for  each $\vp\in(0,1)$,
$$y^\lbb_\vp\to y_\vp\mbox{ in }C([0,T];H^{-1}).$$
\end{lemma}

\pf It suffices to prove the convergence for the solution $z^\lbb_\vp$ to equation \eqref{e3.2} with $\beta$ replaced by $\beta_\lbb$. If we subtract the corresponding equation, we~get
$$\barr{l}
\dd\frac\pp{\pp t}\,(z_\vp-z^\lbb_\vp)+(\vp-\Delta)
((\beta(z_\vp)-\beta_\lbb(z^\lbb_\vp))+\vp(z_\vp-z_\vp^\lbb))
\vspp
\qquad  -{\rm div}(D(z_\vp-z^\lbb_\vp))
+\dd\frac12\,(\mu-\vp^2-(W_\vp)_t)(z_\vp-z^\lbb_\vp)=0,\vspp
(z_\vp-z^\lbb_\vp)(0)=0.\earr$$
Applying $\<z_\vp-z^\lbb_\vp,\cdot\>_{-1,\vp}$ to this equation   and integrating on $(0,t)$, we get
$$\barr{l}
\dd\|(z_\vp-z^\lbb_\vp)(t)\|^2_{-1,\vp}+
\int^t_0\int_{\rr^d}(\beta(z_\vp)-\beta_\lbb(z^\lbb_\vp)
+\vp(z_\vp-z^\lbb_\vp))
(z_\vp-z^\lbb_\vp)ds\,d\xi\vspp
\le C_\vp\dd\int^t_0\|z_\vp(s)-z^\lbb_\vp(s)\|^2_{-1,\vp}ds
+\dd\int^t_0\<{\rm div}\,D(z_\vp-z^\lbb_\vp),z_\vp-z^\lbb_\vp\>_{-1,\vp}ds\vspp \le C_\vp\dd\int^t_0\|z_\vp(s)-z^\lbb_\vp(s)\|^2_{-1,\vp}ds
+C^1_\vp\dd\int^t_0|z_\vp(s)-z^\lbb_\vp(s)|_2
\|z_\vp(s)-z^\lbb_\vp(s)\|_{-1,\vp}ds.\earr$$
This yields
$$\barr{ll}
\|z_\vp(t)-z^\lbb_\vp(t)\|^2_{-1,\vp}
\vsp
\le C^\vp_2\(\dd\int^t_0\|z_\vp(s)-
z^\lbb_\vp(s)\|^2_{-1,\vp}ds
+\dd\int^t_0\int_{\rr^d}|\b(z_\vp)-
\b_\lbb(z^\lbb_\vp)|^2d\xi ds\)  \earr$$
Taking into account that, as easily seen for each $\vp\in(0,1)$, $\{z^\lbb_\vp\}$ is bounded in $L^2((0,T)\times\rr^d)$ and $\beta_\lbb(z_\vp)\to\beta(z_\vp)$, a.e. in $(0,T)\times\rr^d$ as $\lbb\to0$, and     $|\beta_\lbb(z^\lbb_\vp)|\le|\beta(z^\lbb_\vp)|
\le K(1+|z^\lbb_\vp|)$, we infer by  Lebesgue's dominated convergence theorem that, for $\lbb\to0$,
$$z^\lbb_\vp(t)\to z_\vp(t)\mbox{ in }H^{-1}\ \mbox{ uniformly on }[0,T],$$as claimed.\hf

\begin{lemma}\label{l3.5} Let $x\in D(A)\cap D(A_1)$.     Then $y_\vp\in L^\9((0,T)\times\rr^d)\cap L^\9(0,T;L^1)$ and\vspace*{-3mm}
\begin{eqnarray}
&\label{e3.25}
\dd\sup_{\vp\in(0,1)}\{|y_\vp|_{L^\9((0,T)
\times\rr^d)}\}\le C(1+|x|_\9),\\[1mm]
&\sup_{\vp\in(0,1)}|y_\vp|_{L^\9(0,T;L^1)}\le C(|x|_1+1),
\label{e3.24a}\end{eqnarray}
where $C$ is independent of $x$.
\end{lemma}

\pf Let $M= |x|_{\infty}+1$ and $\alpha\in
  C^{1}[0,T]$ be such that $\alpha(0)=0,$ $ \alpha'\geq 0.$ Since $y_\vp $ is a strong solution of \eqref{e3.1} in $H^{-1}$, we have
    \begin{equation}\label{e3.26}
   \barr{l}
    \dd\frac\pp{\pp t}\,(y_\vp-M-\alpha(t))
    -e^{-W_\vp}
    \Delta\(\beta(e^{W_\vp}y_\vp)+\vp e^{W_\vp}y_\vp\)\vsp
    \qquad+\, e^{-W_\vp}\Delta (\beta(e^{W_\vp}(M+\alpha(t)))
    +\vp e^{W_\vp}(M+\alpha(t))) \vsp
     \qquad+\,\vp e^{-W_\vp}
    (\beta(e^{W_\vp}y_\vp)-\beta(e^{W_\vp}(M+\alpha(t))))
    \vsp
     \qquad-\,e^{-W_\vp}{\rm div}\(e^{W_\vp}D(y_\vp-M-\alpha(t))\)\vsp \qquad+\,\dd\frac12\,\mu(y_\vp-M-\alpha(t))=  F_\vp-\alpha',\earr \end{equation}
    where
    \begin{equation}\label{e3.27}
   \barr{ll}
    F_\vp=\!\!\!& e^{-W_\vp}{\rm div}(De^{W_\vp}(M+\alpha(t)))
  -\vp e^{-W_\vp}\beta(e^{W_\vp}(M+\alpha(t)))\vsp
    &     -\dd\frac12\,\mu(M+\alpha(t))
     +e^{-W_\vp}\Delta\beta(e^{W_\vp}(M+\alpha(t)))\vsp
     &
    +\,\vp(M+\alpha(t))e^{-W_\vp} \Delta(e^{W_\vp})
     ,\earr\end{equation}
     and $\a$ will be chosen below, so that
     $$F_\vp-\a'\le0.$$
   To make clear the argument, we shall first prove    \eqref{e3.25} under the condition
   \begin{equation}\label{e3.28}
   \dd\frac{\pp y_\vp}{\pp t}, \beta(e^{W_\vp}y_\vp),\Delta(\beta(e^{W_\vp}y_\vp)+\vp e^{W_\vp} y_\vp)\in L^1((0,T)\times\rr^d).\end{equation}

    Now, we multiply \eqref{e3.26} by ${\rm sign}(y_\vp-M-\alpha(t))^+$ and integrate over  \mbox{$(0,t)\times\rr^d$.}
    We note here that, by \eqref{e1.4}, \eqref{e1.5}, we have that $e^{-W_\vp}\D(\beta(e^{W_\vp}
    (M+\alpha(t))))\in L^1$ and that, after this multiplication, all terms on the left hand side of \eqref{e3.26} become integrable, because of \eqref{e3.28} and, since $\b$ is increasing,   and satisfies \eqref{e1.3}--\eqref{e1.4}. By the monotonicity of $\beta$, and by the elementary inequality
    \begin{equation}\label{e3.29}
    \int_{\rr^d}\Delta z\,{\rm sign}\, (z-M_1)^+d\xi\le0,\ \ff z\in H^1 \mbox{ with } \Delta z\in L^1(\rr^d),\ M_1\ge0.\end{equation}
    we have, because 
    $$\barr{ll}
    {\rm sign}(y_\vp-M- \a(t))^+\!\!\!&=
    {\rm sign}(\beta(e^{W_\vp}y_\vp)-\beta(e^{W_\vp}(M+\a(t))))^+ \vsp&= {\rm sign}((\b+\vp I)(e^{W_\vp}y_\vp)-(\b+\vp I)(e^{W_\vp}(M+\a(t))))^+,\earr$$ where $I(r)=r,$ $r\in\rr$,
    
    \begin{equation}\label{e3.30}
   \barr{lcl}
   J(t):= -\dd \int_{\rr^d} e^{-W_\vp}
    [\Delta(\beta(e^{W_\vp} y_\vp)+\vp e^{W_\vp}y_\vp)\vspp
    -\Delta(\beta(e^{W_\vp}
    (M+\a(t)))+\vp e^{W_\vp}(M+\a(t)))\vspp
    -\vp(\b(e^{W_\vp}y_\vp)-\b(e^{W_\vp}(M+\a(t))))]
    {\rm sign}(y_\vp-M-\a(t))^+ d\xi \vspp
    \ge-
    \dd 2\int_{\rr^d}e^{-W_\vp}\nabla
    [\beta(e^{W_\vp} y_\vp)+\vp e^{W_\vp}y_\vp
    -\beta(e^{W_\vp}(M+\a(t)))\vspp
   -\vp e^{W_\vp}(M+\a(t))]\cdot\nabla W_\vp\,{\rm sign}(y_\vp-M- \a(t))^+  d\xi\vspp
   +\dd\int_{\rr^d}\Delta(e^{-W_\vp})
    (\beta(e^{W_\vp}y_\vp)
    -\beta(e^{W_\vp}(M+\a(t)))\vspp
        +\vp e^{W_\vp}(y_\vp-M-\a(t)))
    {\rm sign}(y_\vp - M-\alpha(t))^+
    d\xi\vsp
    =-\dd\int_{\rr^d}\D(e^{-W_\vp})
    ((\b+\vp I)(e^{W_\vp}y_\vp)
    -(\b+\vp I)(e^{W_\vp}(M+\a(t))))^+d\xi\vsp
    \ge    -(\beta'(e^{-\|W\|_\9}M)+1)e^{\|W\|_\9}
    \|e^{W_\vp}\D(e^{-W_\vp})\|_\9
    %\|e^{W_\vp}\D(e^{-W_\vp})\|_\9
    \dd\int_{\rr^d}(y_\vp-M-\a(t))^+d\xi,
    \earr\hspace*{-8mm}
    \end{equation}
    where, in the last step, we used that on $\{y_\vp-M-\alpha(t)>0\}$ by the mean value theorem and \eqref{e1.5}, we have

    $$\barr{l}
     \beta(e^{W_\vp}y_\vp)-\b(e^{W_\vp}(M+\alpha(t)))
     \le\beta'(e^{W_\vp}(M+\alpha(t)))
     \cdot e^{W_\vp}(y_\vp-M-\alpha(t))\vsp
     \le\beta'(e^{-\|W\|_\9}M)e^{\|W\|_\9}(y_\vp-M-\alpha(t)).\earr$$
      This yields

    \begin{equation}\label{e3.31}
    \barr{r}
    \dd\int^t_0J(s)ds
    \ge -(\beta'(e^{-\|W\|_\9}M)+1)e^{\|W\|_\9}
     (\|\Delta W\|_\9+\|\nabla W\|^2_\9)\\
    \cdot\dd\int^t_0|(y_\vp-(M+\a(s)))^+|_1)ds,\earr
    \end{equation}
     where   $\|\cdot\|_\9$ is the norm of $L^\9((0,T)\times\rr^d).$ (Here and everywhere in the following we shall denote by $C$ several positive constants independent of $W$ and $\vp$.)
     We also have, since $\pp_if\,{\rm sign}\,f^+=\pp_if^+$,

    \begin{equation}\label{e3.32}
 \barr{r}
   \dd\int^t_0\!\!
   \int_{\rr^d}e^{-W_\vp}{\rm div}
   (e^{W_\vp} D(y_\vp{-}M{-}\a(s))){\rm sign}\,
   (y_\vp{-}M{-}\a(s))^+ds\,d\xi\vspp
   =\dd\int^t_0\int_{\rr^d}\nabla W_\vp\cdot D(y_\vp{-}M{-}\a(s))^+  ds\,d\xi .
  \earr\end{equation}
   Taking into account that
   $$\barr{r}
   \dd\int_{\rr^d}
   \frac{\pp }{\pp t}\,(y_\vp(t,\xi)-M-\a(t))
   {\rm sign}(y_\vp(t,\xi)-M-\a(t))^+d\xi\\
   \qquad=\dd\frac d{dt}\,|(y_\vp(t)-M-\a(t))^+|_1,\mbox{ a.e. }t\in(0,T),\earr$$
   after some calculations involving \eqref{e3.26}--\eqref{e3.32}, assuming that $F_\vp\le\a'$, we obtain that
\begin{equation}
  \barr{l}
    |(y_\vp(t)-M-\a(t))^+|_1
    \label{e3.33}
  \le
    \dd\int^t_0\int_{\rr^d}
    ((L+1)(\|\D W_\vp\|_\9\vsp\qquad
   +\|\nabla W_\vp\|^2_\9)
    +\nabla W_\vp\cdot D)(y_\vp-M-\a(s))^+ ds\,d\xi\vsp\qquad
    \le(\beta'(e^{-\|W\|_\9}M)+1)e^{\|W\|_\9}
    (\|\D W_\vp\|_\9+\|\nabla W_\vp\|^2_\9)\vsp\qquad
     +\|\nabla W_\vp\|_\9\|D\|_\9
    \dd\int^t_0
    |(y_\vp(s)-M-\a(s))^+|_1)ds.
    \earr\end{equation}
    By \eqref{e3.33}, it follows that
        \begin{equation}\label{e3.34}
    |(y_\vp(t)-M-\a(t))^+|_1=0\end{equation}
if $F_\vp\le \a'$, a.e. in $(0,T)\times\rr^d.$ To find $\a$ so that this holds, we set
$$\barr{r}
C:=e^{\|W\|}(\|{\rm div}\,D\|_\9+\|D\|^2_\9+\|\mu\|_\9+2
+a_1\vsp
+a_1a_2)(\|\D W\|_\9+\|\nabla W\|^2_\9+1).\earr$$
Then, by assumptions \eqref{e1.3}, \eqref{e1.4}, and  an elementary calculation, we have
   $$    F_\vp \le C (M+\a(t))=\a'(t),$$
    if   $\alpha(t) = M(\exp(C  t)-1)$,  and so \eqref{e3.34} holds. Hence $$y_\vp(t)\le M+\a(t)\le   M+\a(T)<\9,\ \ff t\in[0,T].$$Since the function $r\mapsto -\b(-r)$, $r\in\rr$, enjoys the same properties as $\b$, by a symmetric argument  we get
 $$y_\vp(t)\ge -M-\a(t),\ \ff t\in[0,T],$$ and so \eqref{e3.25} follows.

    To remove condition \eqref{e3.28},
    we are going to approximate \eqref{e3.26}
    by the finite difference scheme \eqref{3.38} below.  To this end, let us first recall that  $A_1$ is the $L^1$-closure of
    $$A_0z=-\Delta(\beta(z)+\vp z)
    +\vp\beta(z)-{\rm div}(Dz)+\frac12\,\mu z,\ z\in D(A_0)$$(see \eqref{e3.13prim}).   Moreover, by \eqref{e3.10secund} for each $f\in  L^1\cap L^\9$ $(\subset H^{-1})$ and $\lbb>\lbb_0$, the equation
   \begin{equation}
   \label{e3.31secund}
   \lbb z+A_0z=f\end{equation}
   has a unique solution $z\in D(A_0)\cap L^\infty\subset L^1\cap H^{1}\cap L^\9$ and $z,\beta(z)\in H^{1}\cap L^1$.
      To see that  indeed  we also have that $z\in L^\9$, we first note that, for all $M\in(0,\9)$, $z\in H^1$, $(z-M)^+=z-z\wedge M\in H^1$ and that it is easy to see that (cf.~\eqref{e3.29})
   \begin{equation}
   \label{e3.31tert}
   \int\Delta(z-M){\rm sign}(z-M)^+d\xi\le0.
   \end{equation}
   Choosing $M=|f|_\9$ and $\lambda\in(0,\9)$ large enough, we have for the solution $z$ of \eqref{e3.31secund} that
   $$\barr{c}
   \lambda(z-M)-\Delta(\beta(z)-\beta(M)+\vp (z-M))
   -{\rm div}(D(z-M))\vsp
   +\dd\frac\mu2\,(z-M)=f-\lambda M+M\,{\rm div}\,D-\dd\frac{\mu}2\,M\le0.\earr$$
   Multiplying by ${\rm sign}(z-M)^+$ and integrating over $\rr^d$ by \eqref{e3.31secund}, it follows that
   $$\lambda\int_{\rr^d}(z-M)^+d\xi+\frac12\int\mu(z=M)^+d\xi\le0,$$hence $z\le M$. Since $r\mapsto-\beta(-r),$ $r\in\rr,$ enjoys the same properties as $\beta$, by symmetry we get $z\ge-M,$ so $z\in L^\9$.
    Hence
    \begin{equation}\label{3.37}
    (\lbb I+A_1)^{-1}(L^1\cap L^\9)\subset D(A_0)\cap L^\infty\subset L^1\cap H^1\cap L^\9,\ \ff\lbb>\lbb_0.\end{equation}

    Now, let us show that the solution $z_\vp $ constructed in Lemma \ref{l3.1} is also the limit of another, for our purpose more convenient finite difference scheme. To this end, define for $h\in(0,1)$ and $0\le i\le N-1$, with $N:=\left[\frac Th\right]$,
    
    $$\nu^h_i:=\frac1h(e^{-W_{i+1}}-e^{-W_i})+(W_\vp)_t(ih)e^{-W_i},$$
    where $W_i:=W_\vp(ih)$. Now, consider the finite difference approximation scheme (again setting $\wt u_i:=\wt u^h_i$)
     \begin{equation}\label{3.36z}
    \barr{l}
    \dd\frac1h(\wt u_{i+1}-\wt u_i)+A_1\wt u_{i+1}+\Lambda(ih)\wt u_{i+1}+\nu^h_i\wt u_{i+1}=0,\vsp
    \wt u^h_0=u_0=x.\earr\end{equation}
    If $u_i:=u^h_i$ is as in \eqref{e5.10prim}, then
    $$\frac1h(u_{i+1}-u_i)+A_1u_{i+1}+\Lambda(ih)u_{i+1}+\nu^h_iu_{i+1}+\eta_i(h)=0,$$
    where $\eta_i(h)=-\nu^h_i\to0$ in $L^1$, uniformly on $[0,T]$ as $h\to0$. Hence, by the same arguments to prove that the schemes \eqref{5.10} and \eqref{e5.10prim} in the proof of Lemma \ref{l5.1} render  the same limit, we obtain that
    $$\lim_{h\to0}\wt u^h=z_\vp \mbox{ in }L^1\mbox{ and $H^{-1}$ uniformly on $[0,T]$}.$$
Setting $y_i:=y^h_i=e^{-W_i}\wt u_i$, we conclude that
\begin{equation}\label{3.36zz}
   \lim_{h\to0} y^h_\vp =y_\vp \mbox{ in }L^1\mbox{ and $H^{-1}$ uniformly on $[0,T]$},\end{equation}
    and, for $0\le i\le N-1$, $N:=\left[\frac Th\right],$
    \begin{equation}\label{3.36zzz}
    \barr{l}
    \dd\frac1h(y_{i+1}-y_i)
    +e^{-W_{i+1}}A_1(e^{W_{i+1}}y_{i+1})=0,\vsp
    	y_0=x_0,\earr\end{equation}
    where $y^h_\vp(t):=y_i$ for $t\in[ih,(i+1)h).$ Since $x\in L^1\cap L^\9$, by \eqref{3.37} we have that $e^{W_i}y_i\in D(A_0)\cap L^\9,$ $0\le i\le N$. So, in \eqref{3.36zzz} we may replace $A_1$ by~$A_0$.

    Now, the approximating scheme \eqref{3.36zzz} can be written as
    \begin{equation}\label{3.38}
    \barr{l}
    \dd\frac1h\,(y_{i+1}-y_i-
    (\alpha(ih)-\alpha((i-1)h) ))
    +e^{-W_{i+1}}(A_0(e^{W_{i+1}}y_{i+1})\\
    \qquad-A(e^{W_{i+1}}(M+\alpha(ih))))
    =F^i_\vp-\dd\frac1h\,
    (\a(ih)-\a((i-1)h))\le0,\earr
    \end{equation}
    where
    $$\barr{ll}
     F^i_\vp=\!\!\!& e^{-W_{i+1}}{\rm div}(De^{W_{i+1}} (M+\alpha(ih)))
    -\vp e^{-W_{i+1}}\beta(e^{W_{i+1}}
    (M+\alpha(ih)))\vsp
  &-\,\dd\frac12\,\mu(M+\alpha(ih))
    +e^{-W_{i+1}}\Delta\beta(e^{W_{i+1}}
    (M+\alpha(ih)))\vsp
    &
    +\,\vp (M+\alpha(ih))e^{-W_{i+1}}
    \Delta e^{W_{i+1}},\earr$$
    where $A(e^{W_{i+1}}(M+\a(ih)))$ is "algebraically" defined as if $A=A_0$, but the argument is not in the domain of $D(A_0)$ (and not even in $D(A_1))$. We note that choosing $\a$ as above, again by \eqref{e1.3}, \eqref{e1.4} and an elementary calculation, we indeed have that the right hand side of \eqref{3.38} is negative.
By \eqref{3.37} we see that $\beta(e^{W_{i+1}}y_{i+1}),
\Delta(\beta(e^{W_{i+1}}y_{i+1})+\vp e^{W_{i+1}}y_{i+1})$ are in  $L^1(\rr^d)$.

Now, we multiply \eqref{3.38} by ${\rm sign}(y_{i+1}-M-\alpha(ih))^+$ and take into account that
\begin{equation}\label{e3.39}
\!\barr{r}
\dd\frac1h\dd\int_{\rr^d}(y_{i+1}-y_i-
(\a(ih)-\a((i-1)h)))
{\rm sign}(y_{i+1}-M-\a(ih))^+d\xi\\
\ge\dd\frac1h\,
(|(y_{i+1}-M-\alpha(ih))^+|_1
-|(y_i-M-\a((i-1)h))^+|_1).\earr\!\!\!\end{equation}
Arguing as in \eqref{e3.30}-\eqref{e3.31},
we get by \eqref{e3.29}

\begin{eqnarray}
&&\hspace*{-7mm}
I^i_1 := -
\dd \int_{\rr^d}e^{-W_{i+1}}[\Delta
(\beta(e^{W_{i+1}}y_{i+1})+\vp e^{W_{i+1}}y_{i+1})
-\Delta(\beta(e^{W_{i+1}}(M+\a(ih)))\nonumber\\[1mm]
&&\hspace*{-7mm}+\vp e^{W_{i+1}}(M+\alpha(ih)))]{\rm sign} (y_{i+1}-(M+\alpha(ih))
 )^+d\xi\label{3.38prim}\\[1mm]
&&\hspace*{-7mm}\ge\dd\int_{\rr^d}
\Delta(e^{-W_{i+1}})
((\b+\vp I)(e^{W_{i+1}}y_{i+1})
-(\b+\vp I)(e^{W_{i+1}}(M+\a(ih))))^+d\xi\nonumber\\[1mm]
&&\hspace*{-7mm}\ge
-(\beta'(e^{-\|W\|_\9}M)+1)e^{\|W\|_\9}
(\|\Delta W\|_\9+\|\nabla W\|^2_\9)
\dd\int_{\rr^d}(y_{i+1}-
M-\alpha(ih))^+d\xi.\nonumber \end{eqnarray}
Similarly, we have
$$
\barr{ll}
I^i_2\!\!\!
&=\dd \int_{\rr^d}
e^{-W_{i+1}}{\rm div}
(e^{W_{i+1}}D(y_{i+1}-M-\alpha(ih)))
{\rm sign} (y_{i+1}-M-\alpha(ih))^+  d\xi\vspp
&=\dd \int_{\rr^d}
(D(y_{i+1}-M-\alpha(ih)))\cdot\nabla W_{i+1}\,{\rm sign}
 (y_{i+1}-M-\alpha(ih))^+ d\xi.\earr$$  This yields
\begin{equation}\label{3.40}
I^i_2\le \|D\|_\9\|\nabla W\|_\9 \int_{\rr^d}
 (y_{i+1}-M-\alpha(ih))^+ d\xi .\end{equation}

Combining   estimates \eqref{3.38}, \eqref{e3.39}, \eqref{3.38prim}, \eqref{3.40} and the facts that \mbox{$\mu\ge0$} and $\b$ is increasing, we get the discrete analogue of \eqref{e3.33}, that~is, for\break $C:=(\beta'(e^{-\|W\|_\9}M)+1)e^{\|W\|_\9}
(\|\D W\|_\9+\|\nabla W\|^2_\9)+\|D\|_\9\|\nabla W\|_\9$,
$$\barr{l}
\dd\frac1h\,(|(y_{i+1}-M-\a(ih))^+|_1
-|(y_1-M-\a((i-1)h))^+|_1)\vspp
\le C
|(y_{i+1}-M-\alpha(ih))^+|_1.\earr$$

Summing up from $i=0$ to $k$, we get
$$\frac 1h\,|(y_{k+1}-M-\a(ik))^+|_1\le C\sum^k_{i=0}|(y_{i+1}-M-\a(ih))^+|_1,$$
which implies, for all $t\in[0,T]$,
$$|(y_\vp^h(t)-M-\a^h(t))^+|_1=0,$$
where $\a^h(t)=ih$ on $[ih,(i+1)h[,$ $0\le i\le N-1.$
Letting $h\to0$ as above, we get \eqref{e3.25}.

To obtain estimate \eqref{e3.24a}, we multiply equation \eqref{3.36zzz}   by sign$\,y_{i+1}$ and integrate over $(0,t)\times\rr^d$.
Then, similarly as above  we find, since $\mu\ge0$ and $\b$ is increasing, that
\begin{equation}\label{3.39prim}
\barr{ll}
\dd\frac1h\,(|y^+_{i+1}|_1-|y^+_i|_1)\!\!\!
&\le\dd\frac1h\int_{\rr^d}(y_{i+1}-y_i)
 {\rm sign}\,y^+_{i+1}d\xi\vsp
 &\le\dd\int_{\rr^d} e^{-W_{i+1}}\D((\b+\vp I)(e^{W_{i+1}}y_{i+1}))
 {\rm sign}\,y^+_{i+1}d\xi\vsp
 &+\dd\int_{\rr^d} e^{-W_{i+1}}{\rm div}(De^{W_{i+1}}y_{i+1}) {\rm sign}\,y^+_{i+1}d\xi\vsp
 &\le\dd\int_{\rr^d}\D e^{-W_{i+1}}(\b+\vp I)
 (e^{W_{i+1}}y_{i+1})
 {\rm sign}\,y^+_{i+1}d\xi\vsp
 &+\dd\int_{\rr^d} \nabla W_{i+1}\cdot Dy^+_{i+1}d\xi\vsp
 &\le C|y^+_{i+1}|_1,\earr
 \end{equation}
 where
 $$C:=(\|\D W\|_\9+\|\nabla W\|^2_\9)(\beta'(e^{-\|W\|_\9}M)+1)
 e^{\|W\|_\9}
 +\|\nabla W\|_\9\|D\|_\9.$$
 Hence, summing from $i=0$ to $k$, we obtain
 $$|y^+_{i+1}|_1\le |x|_1+Ch\sum^k_{i=0}\int_{\rr^d}
  y^+_{i+1}d\xi.$$
 Since $r\mapsto-\b(-r),$ $r\in\rr,$ also fulfills all our assumptions on $\b$, by a symmetry argument we find
 $$|y^-_{i+1}|_1\le |x|_1+Ch\sum^k_{i=0}\int_{\rr^d} y^-_{i+1}d\xi.$$
 This implies that $\ff t\in[0,T]$
 $$|y^h_{\vp}(t)|_1\le|x|_1e^{CT}$$and \eqref{e3.24a} follows letting $h\to0$. \hf

\begin{lemma}\label{l3.4} Let $x\in D(A)\cap D(A_1)$.
    Then there exists an increasing function $C:[0,\9)\to(0,\9)$ such that
     \begin{equation}\label{e3.17a}
\sup_{t\in[0,T]}|y_\vp(t)|^2_2+
\int^T_0\!\!\!\int_{\rr^d}|\nabla\beta(e^{W_\vp}y_\vp)|^2ds\,
d\xi\le C(|x|_\9+|x|_1),\, \ff\vp\in(0,1],\end{equation}for a constant $C>0$, independent of $\vp\in(0,1].$
 \end{lemma}

\pf Clearly, by Lemma \ref{l3.5} we only have to prove the bound in \eqref{e3.17a} for the integral on the left hand side. To this end, we
   multiply \eqref{e3.1} by $\beta(y_\vp)$ and integrate over $(0,t)\times\rr^d$. Taking into account that (see \cite{1}, Lemma 4.4)
    $$\int_{\rr^d}j(y_\vp(t))d\xi
    =\dd\int^t_0{}_{H^{-1}}
    \<\frac{dy_\vp}{ds}\,(s),
    \b(y_\vp(s))\>_{H^1}ds
    +\dd\int_{\rr^d}j(x)d\xi,$$  where $j(r)=\int^r_0\beta(s)ds$, $r\in\rr$, and that
    $$_{H^{-1}}\<\Delta \beta(e^{W_\vp}y_\vp),e^{-W_\vp}
    \beta(y_\vp) \>_{H^1}
    =-\int_{\rr^d}\nabla\beta(e^{W_\vp}y_\vp)\cdot \nabla(e^{-W_\vp}\beta(y_\vp))d\xi,$$
    we get from \eqref{e3.1} that
  \begin{equation}\label{3.45i}
  \barr{l}
   \dd\int_{\rr^d}j(y_\vp(t))d\xi\vsp
   \qquad+
   \dd\int^t_0\int_{\rr^d}[(
   \nabla\beta(e^{W_\vp}y_\vp)
   +\vp\nabla(e^{W_\vp}y_\vp))
    \cdot\nabla(\beta(y_\vp)e^{-W_\vp})]d\xi\,ds   \vspp
    \qquad\le\dd\int_{\rr^d}j(x)d\xi+\int^t_0\int_{\rr^d}
    e^{-W_\vp}{\rm div}(e^{W_\vp}Dy_\vp)
    \beta(y_\vp)d\xi\,ds.
   \earr \end{equation}

  Let us denote the first and second term on the left side of \eqref{3.45i}
 $I_1$ and $I_2$, respectively, and the two on the right $I_3$ and $I_4$. Then
 \begin{equation}\label{3.46i}
 \barr{l}
 I_4=\dd\int^t_0\int_{\rr^d}
y D\cdot(\beta'(y)\nabla y-\b\nabla W)
 d\xi\,ds\vsp
 \qquad\le\|\nabla W\|_\9\|D\|_\9\dd\int^t_0\int_{\rr^d}
 |y_\vp|\,|\beta(y_\vp)|
 d\xi\,ds\vsp
 \qquad+\,\|D\|_\9
 \dd\int^t_0\int_{\rr^d}|y_\vp|
 \beta'(y_\vp)|\nabla y_\vp |
 d\xi\,ds.\earr\hspace*{-16mm}
 \end{equation}
 Obviously, the first integral in
 the preceding line by Lemma \ref{l3.5}
 and \eqref{e1.4} is bounded by $C_1(1+|x|^2_\9)$ with
 a constant
 $C_1>0$ independent of $\vp$.
 Since by \eqref{e1.7}, \eqref{e1.7prim}   we have
 \begin{equation}\label{3.47i}
 \barr{lcl}
 \beta'(y_\vp)&\le&\beta'(e^{-\|W\|_\9+W_\vp}y_\vp)\vsp
 &\le&\vf(e^{-\|W\|_\9})
 \beta'(e^{W_\vp}y_\vp),\earr\end{equation}
 the second integral in the r.h.s. of \eqref{3.46i}, again by Lemma \ref{l3.5} can be bounded by
 $$\delta\int^t_0\int_{\rr^d}|
 \nabla\beta(e^{W_\vp}y_\vp)|^2d\xi\,ds
 +\frac{C_2}\delta\,(1+|x|^2_\9),\ \ff\delta>0,$$
 where $C_2>0$ is a constant independent of $\vp$. So, altogether
 \begin{equation}\label{3.48i}
 I_4\le\delta\int^t_0\int_{\rr^d}|\nabla\b(e^{W_\vp}y_\vp)|^2d\xi\,dx
 +\(\frac{C_2}\delta+C_1\)(1+|x|^2_\9),\ \ff\delta>0.\end{equation}
 Clearly, by \eqref{e1.3},
 \begin{equation}\label{3.49i}
 I_3\le\sup_{r\in|-x|_\9,|x|_\9]}|\b(r)|\
 |x|_1\le C(1+|x|_\9)|x|_1 (<\9).\end{equation}
 Furthermore,
 $$I_2=\int^t_0\int_{\rr^d}\nabla\b(e^{W_\vp}y_\vp)\cdot\nabla(\b(y_\vp)e^{-W_\vp})d\xi\,ds+\vp\wt I_2,$$
 where
 $$\barr{l}
 \hspace*{-6mm}
 \wt I_2:=\dd\int^t_0\int_{\rr^d}\nabla(e^{W_\vp}y_\vp)\cdot
 [\b'(y_\vp)e^{-W_\vp}\nabla y_\vp-\b(y_\vp)e^{-W_\vp}\nabla W_\vp]d\xi\,ds\vsp
 =\dd\int^t_0\int_{\rr^d}
 \nabla(e^{W_\vp}y_\vp)\cdot
 [e^{-2W_\vp}\b'(y_\vp)
 \nabla(e^{W_\vp}y_\vp)\vsp
 -\,e^{-2W_\vp}\b'(y_\vp)y_\vp e^{W_\vp}\nabla W_\vp
 -\b(y_\vp)e^{-W_\vp}\nabla W_\vp]d\xi\,ds\vsp
 \ge-\dd\int^t_0\int_{\rr^d}[|\nabla(e^Wy_\vp)|
 \b'(y_\vp)|y_\vp|+
 |\nabla(e^{W_\vp}y_\vp)|\,
 |\b(y_\vp)|]
 e^{-W_\vp}|\nabla W_\vp|d\xi\,ds.\earr$$
 Since, by Lemma \ref{l3.5}, we have $\sup\limits_{\vp\in(0,1)}\|y_\vp\|_\9<\9$, it follows by \eqref{e1.5} that
 \begin{equation}\label{3.50i}
 \b'(e^{W_\vp}y_\vp)\ge\b'(e^{\|W\|_\9}{\rm sign}\,y_\vp
 \sup\limits_{\vp\in(0,1)}\|y_\vp\|_\9)\ (>0).\end{equation}
 Combining this with \eqref{3.47i}, we conclude that,
 for some increasing functions $\wt C_3,C_3:[0,\9)\to(0,\9)$,
  independent of $\vp$
  \begin{eqnarray}
   \wt I_2
 &\ge&-\dd\int^t_0\int_{\rr^d}
 [\vf(e^{-\|W\|_\9})
 |\nabla\b(e^{W_\vp}y_\vp)|\,|y_\vp|
 \label{3.51i}\\
 &+&(\b'(e^{\|W\|_\9}{\rm sign}\,y_\vp
 \dd\sup_{\vp\in(0,1)}\|y_\vp\|_\9))^{-1}
 |\nabla\b(e^{W_\vp}y_\vp)|\,
 |\b(y_\vp)|]e^{-W_\vp}|\nabla W_\vp|d\xi\,ds\nonumber\\
 &\ge&-\delta\dd\int^t_0\int_{\rr^d}
 |\nabla\b(e^{W_\vp}y_\vp)|^2d\xi\,ds\\
 &-&\dd\frac{\wt C_3(|x|_\9)}\delta
 \dd\int^t_0\int_{\rr^d}
 (|y_\vp|^2+|\b(y_\vp)|^2)
 |\nabla W_\vp|^2d\xi\,ds\nonumber\\
 &\ge&-\delta\dd\int^t_0\int_{\rr^d}
 |\nabla\b(e^{W_\vp}y_\vp)|^2d\xi\,ds
 -\dd\frac{C_3(|x|_\9)}\delta,\ \ff\delta>0,\nonumber\end{eqnarray}
where in the last step we used that, by \eqref{e1.3}, the second integral in the pre\-vious to the last line in \eqref{3.51i} by Lemma \ref{l3.5} is up to a constant, independent of $\vp$, bounded by
$$\int^t_0\int_{\rr^d}|\nabla W|^2d\xi\,ds<\9.$$
Furthermore,
$$\barr{ll}
{\buildrel{\approx}\over{I}}_2:\!\!&
=\dd\int^t_0\int_{\rr^d}
\nabla\b(e^{W_\vp}y_\vp)\cdot\nabla(\b(y_\vp)e^{-W_\vp})d\xi\,ds\vsp
&=\dd\int^t_0\int_{\rr^d}
\nabla\b(e^{W_\vp}y_\vp)\cdot
(e^{-W_\vp}\b'(y_\vp)
\nabla y_\vp-\b(y_\vp)e^{-W_\vp}\nabla W_\vp)d\xi\,ds\vsp
&=\dd\int^t_0\int_{\rr^d}
\b'(e^{W_\vp}y_\vp)
\b'(y_\vp)e^{-2W_\vp}\nabla(e^{W_\vp}y_\vp)
\cdot(\nabla(e^{W_\vp}y_\vp)
-y_\vp e^{W_\vp}\nabla W_\vp )d\xi\,ds\vsp
&-\dd\int^t_0\int_{\rr^d}
\nabla\b(e^{W_\vp}y_\vp)
\b(y_\vp)e^{-W_\vp}\nabla W_\vp\,d\xi\,ds.\earr$$
Since, by \eqref{e1.7}, \eqref{e1.7prim},
\begin{equation}\label{3.52i}
\barr{lcl}
\b'(y_\vp)&\ge&\b'(e^{\|W_\vp\|_\9}y_\vp)
=\b'(e^{\|W_\vp\|_\9-W_\vp}e^{W_\vp}y_\vp)\vsp
&\ge&(\vf(e^{-\|W_\vp\|_\9+W_\vp}))^{-1}
\b'(e^{W_\vp}y_\vp)\vsp
&\ge&\vf(e^{-2\|W\|_\9})^{-1}
\b'(e^{W_\vp}y_\vp),\earr
\end{equation}
we obtain that, for some constant $C_4>0$
and an increasing function $C_5:[0,\9)\to(0,\9)$,
independent of $\vp$,
\begin{equation}\label{3.53i}
\barr{lcl}
{\buildrel{\approx}\over{I}}_2&\ge&C_4
 \dd\int^t_0\int_{\rr^d}
|\nabla\b(e^{W_\vp}y_\vp)|^2
 d\xi\,ds\vsp
 &&-\dd\int^t_0\int_{\rr^d}
|\nabla\b(e^{W_\vp}y_\vp)|
 e^{-W_\vp}|\nabla W_\vp|
 [|y_\vp|\b'(y_\vp)+|\b(y_\vp)|]d\xi\,ds\vsp
 &\ge&(C_4-\delta)
  \dd\int^t_0\int_{\rr^d}
|\nabla\b(e^{W_\vp}y_\vp)|^2
 d\xi\,ds
 -\dd\frac{C_5(|x|_\9)}\delta,\ \ff\delta>0,
\earr\end{equation}
where in the last step we used that
$$e^{-\|W\|_\9}
\dd\sup_{\vp\in(0,1)}
(\|y_\vp\b'(y_\vp)\|_\9+
\|\b(y_\vp)\|_\9)
\dd\int^t_0\int_{\rr^d}
|\nabla W|^2 d\xi\,ds<\9,$$
because of \eqref{e1.3}, \eqref{e1.4} and Lemma \ref{l3.5}. Finally, we note that, by Lemma \ref{l3.5},
\begin{equation}\label{3.54i}
\!\!\!\barr{l}
\dd I_1 \ge{-}\sup\Big\{\b(r)|r\!\in\!\!\Big[{-}
\sup_{\vp\in(0,1)}\!\|y_\vp\|_\9,\!
\sup_{\vp\in(0,1)}\!\|y_\vp\|_\9\Big]\Big\} \!
\dd\sup_{\vp\in(0,1)}\!|y_\vp|_{L^\9(0,T;L^1)}\vsp
\hspace*{5mm}=-C_6(|x|_\9+|x|_1),\earr\!\!\!\end{equation}
for some increasing function $C_6:[0,\9)\to(0,\9)$.

Recalling that $I_2={\buildrel{\approx}\over{I}}_2+\vp\wt I_2$ and combining \eqref{3.48i}, \eqref{3.49i}, \eqref{3.51i}, \eqref{3.53i} and \eqref{3.54i}, the assertion of Lemma \ref{l3.4} is proved.\hf\bk

\n{\bf Proof of existence.} Assume first that   Hypotheses (i)--(iii) hold.
Let $x\in D_0$. It follows, by Lemmas
\ref{l3.5} and \ref{l3.4},
that $\{\beta(e^{W_\vp}y_\vp)\}$
is bounded in $L^2(0,T;H^1)$, $\{y_\vp\}$
is bounded in $L^\9(0,T;L^2)\cap L^\9((0,T)\times\rr^d)$
and $\left\{\frac{dy_\vp}{dt}\right\}$ is bounded
in $L^2(0,T;H^{-1}).$

Moreover, taking into
account that $\nabla\beta(e^{W_\vp}y_\vp)=\beta'(e^{W_\vp}y_\vp)
\nabla(e^{W_\vp}y_\vp)$ and that by assumption
\eqref{e1.5}  and estimate \eqref{e3.25},
    \begin{equation}
    \label{e354prim}
    \beta'(e^{W_\vp}y_\vp)\ge\rho>0
     \mbox{, a.e. in }(0,T)\times\rr^d,\end{equation}
     it follows that $\{y_\vp\}$ is bounded in
     $L^2(0,T;H^1)$. As a matter of fact, we have
     \begin{equation}
     \label{e3.54av}
     \hspace*{-4mm}\barr{l}
     \dd\sup_{\vp\in(0,1]}
 \Bigg\{
 \|y_\vp\|_\9+\|y_\vp\|_{L^\9(0,T;L^1)}
 +\|e^{W_\vp}y_\vp\|_{L^2(0,T;H^1)}
     +\|y_\vp\|_{L^2(0,T;H^1)}\\ \qquad+\|\b(e^{W_\vp}y_\vp)\|_{L^2(0,T;H^1)}+\left\|\dd\frac{dy_\vp}{dt}\right\|_{L^2(0,T;H^{-1})}\Bigg\}
     \le C^*(\oo),\ \oo\in\ooo,\earr\hspace*{-4mm}\end{equation}
where $C^*$ is $\calf$--measurable and $0\!<\!C^*(\oo)\!\le\! C e^{\|W\|_\9}(\exp(\|\nabla W\|_\9{+}\|\D W\|_\9{+}1),$ $\ff\oo\in\ooo$.
     Then, by the Aubin compactness theorem
     (see, e.g., \cite{1}, p. 26),
     $\{y_\vp\}$ is compact in each
     $L^2(0,T;L^2(B_R))$ where
     $B_R=\{\xi\in\rr^d;$ $ |\xi|\le R\}.$
     Hence, for fixed $\oo\in\ooo$ along a subsequence,
     again denoted $\{\vp\}$, we have
    \begin{equation}\label{e3.41}
    \barr{rcll}
    y_\vp&\longrightarrow&y&\mbox{strongly in }L^2((0,T);L^2_{\rm loc}(\rr^d))\vspp
    &&&\mbox{weak-star in }L^\9((0,T)\times\rr^d)\vspp
    &&&\mbox{weakly in }L^2(0,T;H^1),\vspp
    \beta(e^{W_\vp}y_\vp)&\longrightarrow&\eta&\mbox{weakly in }L^2(0,T;H^1)\vspp
    \dd\frac{dy_\vp}{dt}&\longrightarrow&\dd\frac{dy}{dt}&\mbox{weakly in }L^2(0,T;H^{-1})\vspp
    W_\vp&\longrightarrow&W&\mbox{in }C([0,T]\times\rr^d),\earr
    \end{equation}
    and so, letting $\vp\to0$ in equation \eqref{e3.1}, we see that
    \begin{equation}\label{e3.42}
    \barr{l}
    \dd\frac{dy}{dt}-e^{-W}{\rm div}(De^Wy)-e^{-W}\Delta\eta+\frac12\,\mu y=0\mbox{ in }(0,T)\times\rr^d,\vspp
    y(0)=x\mbox{\ \ on }\rr^d.\earr\end{equation}
          Here, we recall that $\D:H^1\to H\1$ is continuous and, by \eqref{e3.6prim}, $u\mapsto {\rm div}(Du)$ is continuous from $L^2$ to $H\1$,   while $\frac d{dt}$ is considered in sense of $H^{-1}$-valued distributions on $(0,T)$ or, equi\-va\-lently, a.e. on $(0,T)$.  Clearly, estimate \eqref{e3.54av}   remains  true for $y$.

    To show that $y$ is a solution to \eqref{e1.13},  it remains to be proven that $\eta=\b(e^Wy)$, a.e. in $(0,T)\times\rr^d$. Since the map $z\to\beta(z)$ is maximal monotone in each $L^2((0,T)\times B_R)$, it is closed and so the latter follows by \eqref{e3.41}. Moreover, if  the solution $y$ to \eqref{e1.13} is unique (we shall see later on that this happens if $\b$ is locally Lipschitz),  it follows that   the sequence $\{y_\vp\}$  arising in \eqref{e3.41}  is independent of $\oo\in\ooo$, and so $y$ is $(\calf_t)_{t\ge0}$--adapted.

    \bk\n{\bf Uniqueness.} Assume   that, besides assumptions (i)-(iii), $\b$ is locally Lipschitz on $\rr$ and consider  $y_1,y_2$ to be two solutions to equation \eqref{e1.13} satisfying \eqref{e2.1}--\eqref{e2.3} and let $z=y_1-y_2$. We~have
  \begin{equation}\label{e3.46}
    \barr{l}
    \dd\frac{\pp z}{\pp t}-e^{-W}{\rm div}(e^WDz)-e^{-W}\Delta
    (\beta(e^Wy_1)-\beta(e^Wy_2))+\frac12\,\mu z=0\vspp\hfill\mbox{ in }(0,T)\times\rr^d,\\
    z(0)=0\mbox{\ \ in }\rr^d.\earr\end{equation}
    Equivalently,
    \begin{equation}\label{e3.47}
    \barr{ll}
    \dd\frac{\pp z}{\pp t}+(I- \Delta)
    (z\eta)\!\!&
    =e^{-W}\,{\rm div}(e^WDz) -e^W\Delta(e^{-W})z\eta\\
 &-\,2\nabla(e^{-W})
    \cdot\nabla(e^Wz\eta)-\dd\frac12\,\mu z+z\eta,\earr\end{equation}
where
$$\eta=\left\{\barr{cl}
\dd\frac{\beta(e^Wy_1)-\beta(e^Wy_2)}{e^Wz}&\mbox{\ \ on }[(\xi,t);z(t,\xi)\ne0],\\
0&\mbox{\ \ on }[(\xi,t);z(t,\xi)=0].\earr\right.$$
We note that, by Hypothesis (ii) \eqref{e1.5},
we have, for some $\a_i=C^i(|x|_1+|x|_\9),$  $i=0,1$,  where $C^i:[0,\9)\to(0,\9)$ are increasing continuous functions,
\begin{equation}\label{e3.48}
0<\alpha_0\le\eta\le\a_1,\mbox{ a.e. in }(0,T)\times\rr^d,\end{equation}
because $\b$ is locally Lipschitz  on $\rr$.

We have $z\in L^2(0,T;H^1(\rr^d))$ and $\frac{\pp z}{\pp t}\in L^2(0,T;H^{-1})$.
We multiply \eqref{e3.47}
by $(I-\Delta)^{-1}z$ and integrate over $\rr^d$ to get

\begin{equation}\label{e3.49}
\barr{l}
 \dd\frac12\,|z(t)|^2_{-1}+
        \dd\int^t_0\int_{\rr^d}
         \eta z^2ds\,d\xi\vspp
         =\dd\frac12\,|z(0)|^2_{-1}+
         \dd\int^t_0\int_{\rr^d}e^{-W}
         {\rm div}(e^W D z)
         (I-\Delta)^{-1}zds\,d\xi\vspp
         -\dd\int^t_0\int_{\rr^d}e^W\Delta(e^{-W})z\eta(I-\Delta)^{-1}zds\,d\xi\vspp
         -2\dd\int^t_0\int_{\rr^d}\nabla(e^{-W})
         \cdot\nabla(e^W z\eta)(I-\Delta)^{-1}zds\,d\xi\vspp
         +\dd\int^t_0\int_{\rr^d}\(-\dd\frac12\,\mu+\eta\)z(I-\Delta)^{-1}zds\,d\xi\vspp
         =\dd\frac12\,|z(0)|^2_{-1}+\int^t_0
         (I_1+I_2+I_3+I_4)ds.\earr\end{equation}
 By the right hand side of \eqref{e3.48}, we get the following estimates
        
         $$\barr{ll}
         |I_1|\!\!
         &\le C|z|_2|z|_{-1}, \vsp
         |I_1|\!\!&\le C\dd\int_{\rr^d}|\b(e^W y_1)-\b(e^Wy_2)|
         \,|(I-\D)^{-1}z|d\xi\vsp
         &\le C|\b(e^Wy_1)-\b(e^Wy_2)|_2|z|_{-1}\le\wt C|z|_2|z|_{-1}\vsp
         |I_3|\!\!&\le\left|2\dd\int_{\rr^d}
         (\b(e^Wy_1)-\b(e^Wy_2)){\rm div}(\nabla(e^{-W})(I-\D)^{-1}z)d\xi\right|\vsp
         &\le C|\b(e^Wy_1)-\b(e^Wy_2)|_{2}|z|_{-1}\le\wt C|z|_2|z|_{-1}\vsp
         |I_4|\!\!&\le C(|z|^2_{-1}+|z|_1|z|_{-1}). \earr$$
         We note also that, by \eqref{e3.48}, we have
         $$|z|_2|z|_{-1}
         \le|\sqrt{\eta}\,z|_2(\sqrt{\a_0})^{-1}|z|_{-1}
         \le\frac12\,|\sqrt{\eta}\,z|^2_2+
         \frac12\,\a_0^{-1}|z|^2_{-1}.$$
         Then, by \eqref{e3.49}, we obtain that
         $$\frac d{dt}\,|z(t)|^2_{-1}\le C_2|z(t)|^2_{-1},\mbox{ a.e. }t>0,$$
         which implies $z\equiv0$, as claimed.

Note also that, by  \eqref{e3.46} and \eqref{e3.49} it follows also that
there exist increasing $C_1,C_2:[0,\9)\to(0,\9)$ such that,
for all $x,\bar x\in x\in D_0$, one~has
\begin{eqnarray}
|y(t,x)-y(t,\bar x)|_{-1} \le  C_1
(|x|_{L^1\cap L^\9}+
|\bar x|_{L^1\cap L^\9})
|x-\bar x|_{-1},\ \ff t\in[0,T],\ \ \ \label{e3.51}\\
|y(t,x)-y(t,\bar x)|_1 \le C_2
(|x|_{L^1\cap L^\9}+
|\bar x|_{L^1\cap L^\9})
 |x-\bar x|_1,\ \ff t\in[0,T].\ \ \ \label{e3.62a}\end{eqnarray}
 Indeed, if one applies \eqref{e3.49} for $z(t)=y(t,x)-y(t,\bar x)$ and uses the above estimates on $I_i$, $i=1,2,3,4,$ and \eqref{e3.48}, one gets   \eqref{e3.51}. To get \eqref{e3.62a}, we multiply \eqref{e3.47} by ${\rm sgn}\,z$ (or, more exactly, by $\calx_\delta(\tau)$, where $\calx_\delta$ is given by \eqref{e3.16az}) and integrate over $\rr^d$.

By \eqref{e3.25}, \eqref{e3.24a} and Lemma \ref{l3.4}, we have
\begin{equation}\label{e3.50a}
\hspace*{-3mm}\barr{r}
|y(t)|_\9+|y(t)|_1+\dd\int^T_0
\int_{\rr^d}|\nabla\beta(y(t,x)(\xi))|^2dt\,d\xi
 \le C_3(|x|_\9+|x|_1),\\ \ff t\in[0,T],\earr
\end{equation}for some increasing functions
 $C_i:[0,\9)\to(0,\9)$, $i=1,2,3$.
This means that, by Lemma \ref{l20},
for all $x\in L^1\cap L^\9$,
 $y=y(t,x)$ extends by density to a strong
 solution to \eqref{e1.13}.  The map $L^1\cap L^\9\ni x\mapsto y(t,x)$ is then Lipschitz on balls in $L^1\cap L^\9$.
  Such a function $y$ satisfies equation \eqref{e1.13}, a.e. on $(0,T)\times\rr^d$, and by \eqref{e3.50a} we have
$$\barr{c}
y\in W^{1,2}([0,T];H^{-1})\cap L^\9((0,T)\times\rr^d),\ \ %\vspp
\beta(e^Wy)\in L^2(0,T;H^1).\earr$$
This completes the proof of Theorem \ref{t2.1}.\hf

   \begin{remark}
   	\label{r3.5}{\rm
   By \eqref{e3.51}
and Lemma \ref{l20}, it follows also that,
for $x\in L^1\cap L^\9$,
there is a unique mild (generalized) solution
\mbox{$y\!\in\! L^\9(0,T;L^1)\cap L^\9((0,T){\times}\rr^d)$}
defined as the limit of mild solutions, that is,
$$y=\lim_{n\to\9}y(\cdot,x_n)\mbox{\ \ in }
L(0,T;L^1)$$for $x_n\to x$ in $L^1$,
 where $\{x_n\}\subset D_0$ and is bounded in
 $L^1\cap L^\9$.}\end{remark}

It should be said also that, in the case where $\b$ is not locally Lipschitz, we do not know  whether we have uniqueness. So, the  sequence $\{y_\vp\}$ arising in \eqref{e3.41}  might depend on the fixed $\oo\in\ooo$ and so   we cannot conclude  that the limit $y$ is $(\calf_t)_{t\ge0}$--adapted.

 %%%%%%%%%%%%%%%%%%%%

%%%%%%%%%%%%%%%%%%%%%%%%

\section{The stochastic equation with nonlinear drift}
\setcounter{equation}{0}

We consider here the equation
\begin{equation}\label{e4.1}
\barr{l}
dX-{\rm div}(a(X))dt-\Delta\beta(X)dt=X\,dW\mbox{ in }(0,T)\times\rr^d ,\vspp
X(0,\xi)=x(\xi),\ \xi\in\rr^d,\earr\end{equation}
where $\beta$ and $W$ are as in Section 1,   while $a:\rr \to\rr^d$  satisfies the following assumption
\bit\item[(iv)] $a$ is  Lipschitzian and $a(0)=0.$  \eit
The strong solution $X$ to equation
\eqref{e4.1} is defined as for equation
\eqref{e1.1}.

For simplicity, we shall use the notations $$u_\xi=\nabla u,\ \ u_{\xi\xi}=\Delta u.$$
By transformation \eqref{e1.12}, we reduce the stochastic equation \eqref{e4.1} to the equation (see \eqref{e1.13})
\begin{equation}\label{e4.2}
\barr{l}
\dd\frac{\pp y}{\pp t}-e^{-W}{\rm div}
(a(e^Wy))-e^{-W}(\beta(e^Wy))_{\xi\xi}
+\dd\frac12\,\mu y=0\mbox{ in }(0,T)\times\rr^d,\vspp
y(t,\xi)=x(\xi).\earr\end{equation}
We have

\begin{theorem}\label{t4.1}
If assumptions {\rm(ii), (iii), (iv)} hold and $\b$ is locally Lipschitz,
for each $x\in D_0$, there is a unique strong solution
$y$ to equation \eqref{e4.2} satisfying
\eqref{e2.1}--\eqref{e2.3}.
Moreover, the process $y$ is
$(\calf_t)_{t\ge0}$-adapted
and, if $x\ge0$, a.e. on $\rr^d$,
then $y\ge0$, a.e. on $(0,T)\times\rr^d$,  and the map $D_0\ni x\to y(\cdot,x)$
is Lipschitz from $H^{-1}$ to $C([0,T],H^{-1})$
on balls in $L^1\cap L^\9$ and
 extends to a strong solution
to \eqref{e4.1} satisfying \eqref{e2.1},
\eqref{e2.3}, for all
$x\in L^1\cap L^\9$. \end{theorem}

\pf Since the proof is essentially
the same as that of Theorem \ref{t2.1},
we only sketch it, by emphasizing, however, the points where arise major differences in the argument.\mk

We consider the approximating equation (see \eqref{e3.1})
\begin{equation}\label{e4.3}
\barr{l}
\dd\frac{\pp y_\vp}{\pp t}-e^{-W_\vp} {\rm div} (a(e^{W_\vp}y_\vp))
-e^{-W_\vp} \beta(e^{W_\vp}y_\vp)
-\vp e^{-W_\vp} (e^{W_\vp}y_\vp)_{\xi\xi}\vspp
\qquad+\vp e^{-W_\vp}\beta(e^{W_\vp}y_\vp)+
\dd\frac12\,\mu y_\vp=0\mbox{ in }(0,T)\times\rr^d,\vspace*{3,5mm}\\
y_\vp(0,\xi)=x(\xi),\ \xi\in\rr^d,
\earr\end{equation}
which, by the same argument as that in the proof of Lemma \ref{l3.1}, has a unique solution $y_\vp$ which satisfies \eqref{e3.3}--\eqref{e3.5}.

We note   that Lemmas \ref{l3.1}, \ref{l3.5}   and  \ref{l3.4}   remain  valid in this case too. Indeed, we note that, instead of \eqref{e3.26} and \eqref{e3.27}, we have

\begin{equation}\label{e4.4}
\barr{l}
\dd\frac\pp{\pp t}\,(y_\vp-M-\a(t))-e^{-W_\vp}
(\beta(e^{W_\vp}y_\vp)+\vp e^{W_\vp}y_\vp)\vspp
\qquad- (\beta(e^{W_\vp}(M+\a(t))-\vp
e^{W_\vp}(M+ \a(t))))_{\xi\xi}\vspp
\qquad+\vp e^{-W_\vp}(\beta(e^{W_\vp}y_\vp)-
\beta(e^{\|W_\vp\|}(M+ \a(t))))\vspp
\qquad-e^{-W_\vp}({\rm div} (a(e^{W_\vp}y_\vp)
-a(e^{W_\vp}(M+ \a(t)))))\vspp
\qquad-\dd\frac12\,\mu(y_\vp-M- \a(t))
=F_\vp-\a'(t),\earr\end{equation}
where
$$\barr{c}
F_\vp= e^{-W_\vp}{\rm div}\,
a(e^{W_\vp}(M+ \a(t)))-\dd\frac12\,(M+\a(t))+e^{-W_\vp}
(\b(e^{W_\vp}(M+\a(t))))_{\xi\xi} \vsp
 -\vp e^{-W_\vp}\beta(e^{W_\vp}(M+ \a(t)))
+\vp(M+\a(t))e^{-W_\vp}(e^{W_\vp})_{\xi\xi}
\earr$$
(or it discretized analogue \eqref{3.38}).

In order to treat the term in $a$ arising in \eqref{e4.3}, we note that

$$\barr{l}
-\dd\int^t_0\!\!\!\int_{\rr^d} e^{-W_\vp}{\rm div}(a(e^{W_\vp}y_\vp)
-a(e^{W_\vp}(M+  \a(s)))
{\rm sign}(y_\vp-(M+ \a(s)))^+)ds\,d\xi\vspp
\qquad-\dd\int^t_0\!\!\!\int_{\rr^d}({\rm div}(a(e^{W_\vp}y_\vp)
-a(e^{W_\vp}(M+\a(s)))
e^{-W_\vp}) \vspp
\hfill{\rm sign}(e^{W_\vp}y_\vp
-e^{W_\vp}(M+\a(s))e^{-W_\vp})^+)ds\,d\xi\vspp
\qquad+\dd\int^t_0\!\!\!\int_{\rr^d}(e^{-W_\vp})_\xi
\cdot (a(e^{W_\vp}y_\vp)
-a(e^{W_\vp}(M+\a(s)))\vspp
\hfill
{\rm sign}(e^{W_\vp}y_\vp-e^{W_\vp}(M+\a(s)))^+ds\,d\xi\vspp
\qquad\le L\dd\int^t_0\int_{\rr^d}
|(e^{-W_\vp})_\xi|
(y_\vp-M-\a(s))^+ds\,d\xi,\earr$$
because $a$ is Lipschitz and
\begin{equation}\label{e4.5}
\int_{\rr^d}(e^{-W_\vp} (a(u)-a(v)))_\xi{\rm sign}(u-v)d\xi=0\end{equation}for $u=e^{W_\vp}y_\vp$ and $v=e^{W_\vp}(M+\a(t))$.  To prove \eqref{e4.5}, we consider the   approximation $\calx_\delta$ of the signum function  defined by \eqref{e3.16az}.
We have

$$\barr{ll}
H_\delta(t)\!\!\!
&=\dd\int_{\rr^d}{\rm div}(e^{-W_\vp}(a(u)-a(v)))
\chi_\delta(u-v)d\xi\vspp
&=-\dd\int_{\rr^d}e^{-W_\vp}(a(u)-a(v))\cdot(u-v)_\xi\chi'_\delta(u-v)d\xi\vspp
&=-\dd\frac1\delta\int_{[|u-v|\le\delta]}
e^{-W_\vp}(a(u)-a(v))\cdot(u-v)_\xi d\xi.\earr$$
For $\delta\to0$, we get
$$
\dd\lim_{\delta\to0}H_\delta(t)=\dd\int_{\rr^d}e^{-W_\vp}
{\rm div}(a(u(t,\xi))-a(v(t,\xi))){\rm sign}(u(t,\xi)-v(t,\xi))d\xi$$while
$$
|H_\delta(t)| \le {\rm Lip}(a)\int_{[|u-v|\le\delta]}e^{-W_\vp}|(u-v)_\xi|d\xi. $$
This yields
$$\limsup_{\delta\to0}|H_\delta(t)|\le\int_{[|u-v|=0]}e^{-W_\vp}|(u-v)_\xi|d\xi=0,$$
because $(u-v)_\xi=0$ on $\{\xi;\ (u-v)(\xi)=0\}.$ (We recall that $u,v\in H^1$.)

Then  estimate \eqref{e3.32} with $a$ in place of $D$ remains true in this case.

Multiplying \eqref{e4.4} by sign$(y_\vp-M- \a(t))^+$ and integrating on $(0,t)\times\rr^d$, we get by \eqref{e3.31} an estimate of the form \eqref{e3.33} from which we infer that
$$|(y_\vp(t)-M-\a(t))^+|_1=0,\ t\in(0,T),$$for $\a$   chosen as in the proof of Lemma \ref{l3.5} and so
$$y_\vp\le M+\a(t),\mbox{ a.e. in }(0,T)\times\rr^d,$$
and, similarly,
$$y_\vp\ge -M-\a(t),\mbox{ a.e. in }(0,T)\times\rr^d.$$

Taking into account that
$$\barr{l}
\dd\int^t_0\int_{\rr^d}
e^{-W_\vp}{\rm div}(a(e^{W_\vp}y_\vp))_\xi
\beta(y_\vp)ds\,d\xi%\vspp\qquad
=-\dd\int^t_0\int_{\rr^d}a(e^{W_\vp}y_\vp)
\cdot(e^{-W_\vp}\beta(y_\vp))_\xi ds\,d\xi\vspp
\qquad\le C\dd\int^t_0\int_{\rr^d}(|e^{W_\vp}y_\vp|
(|y_\vp|^m
|(e^{-W_\vp}))_\xi|+e^{-W_\vp}\beta'(y_\vp)|\nabla y_\vp|)ds\,d\xi,\earr$$ and, recalling that $\sup\limits_{\vp>0}\{|y_\vp|_\9\}<\9$, it follows
  as in the proof of Lemma \ref{l3.4}  that estimate \eqref{e3.17a}  holds in this case too. Hence,   there is $y\in C([0,T];L^2_{\rm loc})\cap L^\9((0,T)\times\rr^d)\cap L^2(0,T;H^1)$ such that \eqref{e3.41} holds. Moreover, we have, for $\vp\to0$,
$$a(e^{W_\vp}y_\vp)\to a(e^Wy)\mbox{ in }L^2((0,T);L^2_{\rm loc})$$and so, for $\vp\to0$
$$ {\rm div}(a(e^{W_\vp}y_\vp)) \to {\rm div}
(a(e^Wy))\mbox{ in }L^2([0,T);H^{-1}_{\rm loc}).$$
Then letting $\vp\to0$ in \eqref{e4.3}, we see that $y$ is a solution to equation \eqref{e4.2} satisfying \eqref{e2.1}-\eqref{e2.3}. Moreover, multiplying \eqref{e4.3} by sign$y_\vp$ and taking into account that, as seen earlier,
$$\int_{\rr^d}e^{-W_\vp}{\rm div}(a(e^{W_\vp}y_\vp)) {\rm sign}\,y_\vp d\xi\le C\int_{\rr^d}|e^{W_\vp}y_\vp|d\xi,$$we get as in the proof of Lemma \ref{l3.5} that
$$|y_\vp(t)|_1\le C|x|_1,\ \ff t\in[0,T],$$
where $C$ is independent of $\vp$.

\bk\n{\bf Uniqueness.}  If $\b$ is locally Lipschitz and $y_1,y_2$ are solutions to \eqref{e4.1}, for $z=y_1-y_2$, we get
(see~\eqref{e3.46})
$$\barr{l}
\dd\frac{\pp z}{\pp t}- {\rm div}(a(e^Wy_1)-a(e^Wy_2)) -e^{-W} (\beta(e^Wy_1)-\beta(e^Wy_2))_{\xi\xi}+\frac12\,\mu z=0\vspp
z(0)=0,\earr$$
and, arguing as in the proof of uniqueness in Theorem \ref{t2.1}, we get $z\equiv0$. If $\b\in L^1_{\rm loc}(\rr)$, then, multiplying scalarly in $L^2$
by $(I-\Delta)^{-1}z$  and using the local Lipschitzianity of $\beta$ and $a$,
we get as above  the    estimates \eqref{e3.51}--\eqref{e3.50a}.~\hf\bk

By Theorem \ref{t4.1}, we have

\begin{corollary}\label{c4.2}
If assumptions  {\rm (ii), (iii), (iv)} hold and $\b$ is locally Lipschitz, then
for each $x\in D_0$
there is a unique strong solution
$X$ to the stochastic equation \eqref{e4.1},
which satisfies
\begin{equation}\label{e2.6axx}
X e^{-W}\in W^{1,2}([0,T];H^{-1}),\ \pas,\end{equation}
and  $X\ge0$, a.e. on $(0,T)\times\rr^d\times\ooo$
if $x\ge0$, a.e. on $\rr^d$. Moreover,   the map $x\mapsto X(t,x)$ is $H^{-1}$-Lipschitz from balls in $L^1\cap L^\9$ to $C([0,T];H^{-1})$.
\end{corollary}

\begin{remark}\label{r4.3}
\rm If  $a$ is not Lipschitz, one cannot expect a strong solution for equation \eqref{e4.1}.   In the deterministic case, if $\beta\equiv0$, equation \eqref{e4.1} reduces to a first order quasilinear equation previously studied by S.~Kruzkov \cite{13a} (see, also, \cite{8a}, \cite{8aa}), who introduced and proved existence of a generalized solution involving the so-called "entropy" conditions. (See also \cite{2} for the case where $\b$ is present.) So, also in this case, one might expect to have a generalized solution in sense of Kruzkov, but this remains to be done.    \end{remark}

\section{Appendix}

\setcounter{equation}{0}

Here, we shall briefly review a few definitions and results pertaining the nonlinear Cauchy problem in Banach spaces for quasi-$m$-accretive operators.

Let $X$ be a Banach space with the norm denoted $\|\cdot\|_X$. A nonlinear ope\-rator $A:D(A)\subset X\to X$ (possibly multivalued) is said to be accretive~if
$$\|x_1-x_2+\lbb(y_1-y_2)\|_X\ge\|x_1-x_2\|_X,\ \ff\lbb>0,\ \ff y_i\in Ax_i,\ i=1,2,$$and quasi-accretive if $A+\a I$ is accretive for some $\a>0$. Equivalently,
$$_X(y_1-y_2,\eta)_{X'}\ge0,\mbox{ for some }\eta\in J(x_1-x_2),$$where $J:X\to X'$ is the duality map of the space $X$. (Here, $X'$ is the dual of $X$.) The operator $A$ is said to be {\it m-accretive} if the range $R(\lbb I+A)$ of $\lbb I+A$ is all of $X$ for all $\lbb>0$ and quasi $m$-accretive if $R(\lbb I+A)=X$ for $\lbb>\lbb_0>0$.

If $A$ is quasi $m$-accretive, $u_0\in\ov{D(A)}$ and $g\in C([0,T];X)$, then the Cauchy problem
\begin{equation}\label{e5.1}
\barr{l}
\dd\frac{du}{dt}+Au\ni g\mbox{ in (0,T)},\vspp
u(0)=u_0,\earr\end{equation}
has a unique mild solution $u\in C([0,T];X)$ defined by
\begin{equation}\label{e5.2}
u(t)=\lim_{h\to0} u^h(t)\mbox{ strongly in $X$ and uniformly on }[0,T],\end{equation}
\begin{equation}\label{e5.3}
\barr{l}
u^h(t)=u^h_i\mbox{ for }t\in[ih,(i+1)h],\vspp
\dd\frac1h\,(u^h_{i+1}-u^h_i)+Au_{i+1}\ni\dd\frac1h
\int^{(i+1)h}_{ih}g(t)dt,\vsp
 \hfill i=0,1,...,N-1,\mbox{ with }N=\left[\frac Th\right],\vspp
u^h_0=u_0.\earr\end{equation}
(See, e.g., \cite{1}, Section 4.1, Corollary \ref{c4.2}.) (For $g\equiv0$, this is just the Crandall-Liggett exponential formula.) Moreover, if the space $X$ is reflexive and $g\in W^{1,1}([0,T];X)$, then $u$ is a strong absolutely continuous solution to \eqref{e5.1}, that is, it satisfies a.e. \eqref{e5.1} and
\begin{equation}\label{e5.4}
u\in W^{1,\9}([0,T];X) ,\ Au\in L^\9(0,T;X).\end{equation}
Finally, if $X$ is uniformly convex, then $\frac d{dt}\,u(t)$ is continuous from the right.

We consider now the Cauchy problem
\begin{equation}\label{e5.5}
\barr{l}
\dd\frac{du}{dt}\,(t)+Au(t)+\Lambda (t)u(t)=0,\ \ff t\in(0,T),\vspp
u(0)=u_0,\earr\end{equation}
where $A$ is  quasi-$m$-accretive, $u_0\in\ov{D(A)}$ and $\Lambda \in C([0,T];L(X,X))$. Since it is enough for the applications in this paper, let us for simplicity assume that $A$ is single-valued. We have

\begin{lemma}\label{l5.1} The Cauchy problem
\eqref{e5.5} has a unique mild solution
$u\in C([0,T];X)$ and $u$ is given as
the limit in \eqref{5.9} of the finite
difference scheme \eqref{e5.10prim} below. Moreover, if $u_0\in D(A)$ and
\begin{equation}\label{e5.6}
\|\Lambda (t)-\Lambda (s)\|_{L(X,X)}\le L|t-s|,\ \ff s,t\in[0,T],\end{equation}
then $u:[0,T]\to X$ is Lipschitz.
\end{lemma}

\pf Consider the operator $\cala:D(\cala)\subset L^1(0,T;X)\to L^1(0,T;X)$ defined by
 $\cala u=g$ if  $u\in C([0,T];X)$
is the mild solution to \eqref{e5.1}. By the existence theory for \eqref{e5.1}, it follows that $R(\lbb I+\cala)=L^1(0,T;X)$, $\ff \lbb>0$, and by \eqref{e5.3} we see that $\cala$ is quasi-accretive. Indeed, if $\lbb_0\ge0$ such that $A+\lbb_0 I$ is $m$-accretive, then by \cite{1}, Theorem 4.1 and Proposition 3.7(iv),  we have for solutions $u,\bar u$ for \eqref{e5.1} with $g,\bar g$, respectively, on the right hand side
$$\|u(t)-\bar u(t)\|_X\le\int^t_0
e^{\lbb_0(t-s)}\|g(s)-\bar g(s)\|_Xds,\ \ff\lbb>0,\ g,\bar g\in L^1(0,T;X),$$which yields
$$\|u-\bar u\|_{L^1(0,T;X)}\le\frac{e^{\lbb_0T}}{\lbb_0}\,\|g-\bar g\|_{L^1(0,T;X)}.$$
Hence $\cala$ is quasi-$m$-accretive.

The operator $\wt \Lambda:L^1(0,T;X)\to L^1(0,T;X)$ defined~by
\begin{equation}\label{e5.7}
(\wt \Lambda u)(t)=\Lambda(t)u(t),\ t\in[0,T],\end{equation}
is linear continuous and this implies that $\cala+\wt \Lambda$ is quasi $m$-accretive in $L^1(0,T;X)$. Hence
there is $\lbb_0>0$ such that $R(\lbb I+\cala+\wt \Lambda)=L^1(0,T;X)$ for $\lbb>\lbb_0>0$.

This means that, for every $g\in C([0,T],X)$, the equation
\begin{equation}\label{e5.7prim}
\barr{l}
\dd\frac{du}{dt}+Au+\lbb u=g(t)-\Lambda(t)u,\ t\in(0,T),\vspp
u(0)=u_0\earr\end{equation}
has a unique mild solution for $\lbb>\lbb_0.$

Now, let us show that this implies that also \eqref{e5.5} has a unique mild solution. This is well known, but we include the proof for the reader's convenience. So, fix $\lbb>\lbb_0$ and let $u,\bar u$ be the unique mild solutions of \eqref{e5.7prim} with $\lbb g$ and $\lbb\bar g$ replacing $g$ on its right hand side, where $g,\bar g\in\chi:=C([0,T];X)$, equipped with the norm $\|\cdot\|_\chi:=\|\cdot\|_{\chi,T}$, where for $t\in[0,T]$
$$\|g\|_{\chi,t}:=\sup\{e^{-\a s}\|g(s)\|_X;\ s\in[0,t]\}$$and $\a>0$ will be chosen later. Then, by \cite{1}, Theorem 4.1 and Proposition  3.7(iv), for all $t\in[0,T]$, it follows that
$$\|u-\bar u\|_{\chi,t}\le\dd\int^t_0 e^{-(\lbb-\lbb_0+\a)(t-s)}
(\lbb\|g-\bar g\|_{\chi,s}ds
+C \|u-\bar u\|_{\chi,s})ds,$$
where $C:=\sup\limits_{t\in[0,T]}\|\Lambda(t)\|_{L(X,X)}.$ Hence, by Gronwall's lemma,
$$\|u-\bar u\|_\chi\le\frac{\lbb e^{CT}}{\lbb-\lbb_0+\a}\,\|g-\bar g\|_\chi.$$
Now, choosing $\a$ large enough, it follows that the map which maps $g$ to the solution $u$ of \eqref{e5.7prim} with $\lbb g$ replacing $g$ on its right hand side, is a strict contraction on $\chi$. Hence, by Banach's fixed point theorem, \eqref{e5.5} has a unique mild solution, $u$.

Moreover, as a mild solution to \eqref{e5.5}, by \eqref{e5.2}and \eqref{e5.3} where $g(t)=\Lambda(t)u(t),$ $u$ satisfies
\begin{equation}
\label{5.9}
u=\lim_{h\to0} u^h(t)\mbox{ strongly in $X$ and uniformly on $[0,T]$},
\end{equation}
where, for $h>0$,
\begin{equation}
\label{5.10}\barr{l}
u^h(t)=u^h_i\mbox{ for }t\in[ih,(i+1)h),\vsp
\dd\frac1h\,(u^h_{i+1}-u^h_i)+Au^h_{i+1}+\frac1h\int^{(i+1)h}_{ih}\Lambda(t)u(t)dt=0,\vsp
\hfill i=0,1,...,N-1,\mbox{ with $N=\left[\frac Th\right]$},\vsp
u^h_0=u_0.\earr
\end{equation}

\noindent As easily seen, we may replace \eqref{5.10} by
\begin{equation}\label{e5.10prim}
\frac1h\,(u^h_{i+1}-u^h_i)+Au^h_{i+1}+\Lambda(ih)u^h_{i+1}=0.
\end{equation}
Indeed, setting $u_i:=u^h_i$, we may rewrite \eqref{5.10} as
\begin{equation}\label{e5.11prim}
\frac{1}{h}\,(u_{i+1} - u_i) + Au_{i+1} +\Lambda(ih)u_{i+1} + \eta_i(h)=0,\end{equation}
where $\|\eta _i\| \leq \delta(h) ,\ \forall i,$ and
$\delta(h)\rightarrow 0$ uniformly on $[0,T]$ as $h= \frac TN$ goes to zero.

Now, if $v= v_i ,\ i=0,1, ... N-1,$ is the solution to \eqref{e5.10prim},
subtracting the  equation  \eqref{e5.10prim} from \eqref{e5.11prim}, we get for $y_i = u_i - v_i$  the equation
$$ y_{i+1}+ h(Au_{i+1} - Av_{i+1}) +h \Lambda(ih) y_{i+1}= y_i - h\eta_i(h)$$
and, by the quasi-accretivity of $A$, this yields
$$\|y_{i+1}\| \leq \mu h \|y_{i+1}\| + \|y_i\| + h\delta(h),\ \forall
i =0,1, .. ., n-1,$$
where $\mu = \lambda + \sup\limits_{t\in[0,T]} \|\Lambda(t)\|_{L(X;X)}.$
This yields for   small enough $h$
$$ \|y_{i+1}\| \leq(1-\mu h)^{-1}(\| y_i\| + h\delta(h))$$
and, taking into account that $y_0=0$ and that  $h=\frac TN,$ we get that for $h$ small enough
$$\|y_{i+1} \|\leq  h\delta(h)(1-\mu h)^{-1} \sum_{1\leq j\leq i}
(1-\mu h)^{-j} \leq \frac{\delta(h)}\mu\(1-\frac{T\mu}N\)^{-N}.$$
Hence $y_i= y^h_i$ goes to zero in $X$ as $h$ goes to zero and
this completes the proof of the equivalence of \eqref{e5.10prim}  and \eqref{5.10}.

Now, we shall prove that, if $u_0\in D(A)$ and \eqref{e5.6} holds, then $u$ is Lipschitz.

By \eqref{e5.6}, we have
\begin{equation}\label{e5.9}
\barr{l}
\|\Lambda(t)u(t)-\Lambda(s)u(s)\|_X\vsp
\qquad\quad
\le L|t-s|\|u\|_{C([0,T];X)}
+\|\Lambda(t)\|_{L(X,X)}\|u(t)-u(s)\|_X\vsp
\qquad\quad\le C_1(|t-s|+\|u(t)-u(s)\|_X),\ \ff s,t\in[0,T].\earr
\end{equation}
We consider now the equation
\begin{equation}\label{e5.11}
\barr{l}
\dd\frac{du_\lbb}{dt}+A_\lbb u_\lbb+\Lambda(t)u=0,\ t\in[0,T],\vsp
u_\lbb(0)=u_0,\earr
\end{equation}
where $A_\lbb=\lbb^{-1}(I-(I+\lbb A)^{-1})$ is the Yosida approximation of $A$. The Cauchy problem has a unique differentiable solution $u_\lbb:[0,T]\to X$ and, since $A_\lbb$ is $\wt\lbb_0$-accretive for some $\wt\lbb_0>0$, we have by \eqref{e5.11}
$$\barr{l}
\dd\frac12\,\frac d{dt}\,\|u_\lbb(t+h)-u_\lbb(t)\|^2_X
\le\|\Lambda(t+h)u(t+h)-\Lambda(t)u(t)\|_X\vsp
\|u_\lbb(t+h)-u_\lbb(t)\|_X+\wt\lbb_0\|u_\lbb(t+h)-u_\lbb(t)\|^2_X,\ t,t+h\in[0,T].\earr$$
By \eqref{e5.9}, this yields
\begin{equation}\label{e5.12}
\barr{ll}
\|u_\lbb(t+h)-u_\lbb(t)\|_X
\le e^{(\wt\lbb_0+C_1)t}\|u_\lbb(h)-u_\lbb(0)\|_X\vsp
\qquad+C\dd\int^t_0e^{(\wt\lbb_0+C_1)(t-s)}(h+\|u(s+h)-u(s)\|_X)ds.\earr
\end{equation}
On the other hand, by \eqref{e5.11} we have
$$\barr{ll}
\dd\frac12\,\frac d{dt}\,\|u_\lbb(t)-u_0\|^2_X
\!\!&\le\wt\lbb_0\|u_\lbb(t)-u_0\|^2_X
+\|A_\lbb u_0\|_X
\|u_\lbb(t)-u_0\|_X\vsp
&+\|\Lbb(t)u(t)\|_X\|u_\lbb(t)-u_0\|_X,\ \ff t\in[0,T].\earr$$
Hence
$$\barr{ll}
\|u_\lbb(t)-u_0\|_X
\!\!&
\le\dd\int^t_0e^{\wt\lbb_0(t-s)}
(\|A_\lbb u_0\|_X+\|\Lambda(s)u(s)\|_X)ds\vsp
&\le C_2 (\|Au_0\|_X+1),\ \ff t\in[0,T].\earr$$
Substituting into \eqref{e5.12}, yields
\begin{equation}\label{e5.13}
\barr{l}
\|u_\lbb(t+h)-u_\lbb(t)\|_X\vsp
\qquad\quad
\le C_3\(h+\dd\int^t_0 e^{(\wt\lbb_0+C_1)(t-s)}(h+\|u(s+h)-u(s)\|_X)\)ds,\vsp
\qquad\quad\ff\lbb>0,\ t,t+h\in[0,T].\earr\end{equation}
On the other hand, since for each $\vp>0$
$$\lim_{\lbb\to0}(I+\vp A_\lbb)^{-1}x=(I+\vp A)^{-1}x,\ \ff x\in H,$$
by the Trotter-Kato theorem for nonlinear semigroups of contractions, we have (see \cite{1}, Corollary 4.5)
$$u_\lbb\longrightarrow v\mbox{\ \ in }C([0,T];X)\mbox{\ \ as }\lbb\to0,$$where $v$ is the solution to
$$\barr{l}
\dd\frac{dv}{dt}+Av+\Lambda(t)u=0,\vsp
v(0)=u_0.\earr$$
By the quasi-accretivity of $A$, it follows that $v=u$.
Then, letting $\lbb\to0$ in \eqref{e5.13}, we get
$$\|u(t+h)-u(t)\|_X\le C_3\(h+\dd\int^t_0(h+\|u(s+h)-u(s)\|_X)ds\)$$and by Gronwall's inequality, we get
$$\|u(t+h)-u(t)\|_X
\le C_4h,\ \ff t,t+h\in[0,T],$$as claimed. This completes the proof.\hf\bk

If the space $X$ is reflexive, we infer that, under the conditions of Lemma \ref{l5.1}, $u\in W^{1,\9}([0,T];X)$ is a.e. differentiable, and satisfies equation \eqref{e5.5}, a.e. on~$(0,T).$   We have, therefore,

\begin{corollary}\label{c5.2} If the space $X$ is reflexive, $u_0\in D(A)$, and $\Lambda$ satisfies \eqref{e5.6}, then the mild solution $u$ to \eqref{e5.5} is a strong absolutely continuous solution, which satisfies~\eqref{e5.4}.
\end{corollary}

It should be mentioned that the latter case applies
to $X=H^{-1}$, but not to $X=L^1$. In the latter case, the solution $u$ is only continuous.

\bk\n{\bf Acknowledgement.} This work was supported by the DFG through CRC 1283. Viorel Barbu was also partially supported by CNCS-UEFISCDI (Roma\-nia), through project PN-III-P4-ID-PCE-2016-0011.

\end{document}